\documentclass[12pt,a4paper]{article}
\usepackage{fullpage,amsfonts,amssymb,latexsym,bbm,amsmath}
\usepackage{graphicx} 
\usepackage{indentfirst}
\usepackage{pdfsync}
\usepackage{stackrel}
\usepackage{xcolor}
\usepackage{graphicx}
\usepackage[all]{xy}
\usepackage{amscd}
\usepackage {stackrel}
\usepackage{ifpdf}
\usepackage{appendix}
\ifpdf%
  \usepackage{pdftricks}
  \begin{psinputs}
    \usepackage{pstricks}
  \end{psinputs}
\else
  \usepackage{pstricks}
  
\fi

\newcommand{\AAA}{\mathbb{A}}
\newcommand{\RR}{\mathbb{R}}
\newcommand{\CC}{\mathbb{C}}
\newcommand{\DD}{\mathbb{D}}

\newcommand{\TT}{\mathbb{T}}
\newcommand{\ZZ}{\mathbb{Z}}

\newcommand{\GG}{\mathbb{G}}
\newcommand{\NN}{\mathbb{N}}

\def\ch{{\hbox{ch}}}

\newcommand\hA{\hat{A}}

\newcommand\hphi{\hat{\phi}}

\newcommand{\Cj}{\mathcal{J}}

\newcommand{\dbar}{\overline{\partial}}

\def\Re{\mathop{\rm Re} \nolimits}

\def\beginpf{\noindent {\sl Proof:\,\,}}

\newtheorem{theorem}{Theorem}
\newtheorem{lemma}{Lemma}
\newtheorem{corollary}{Corollary}
\newtheorem{definition}{Definition}
\newtheorem{proposition}{Proposition}
\newtheorem{remark}{Remark}
\newcommand{\boite}{\mbox{} \hfill \mbox{\rule{2mm}{2mm}}}

\def\endpf{\boite\medskip}

\def\beq{\begin{equation}}
\def\eeq{\end{equation}}
 \def\essup{\text{ess sup}}

\begin{document}

\date{}
\title{\bf{Composition operators on generalized Hardy spaces}}
\author{Sam J. Elliott$^{\hbox{\small{a }}}$, Juliette Leblond$^{\hbox{\small{ b}}}$, Elodie Pozzi$^{\hbox{\small{ c}}}$ and Emmanuel Russ$^{\hbox{\small{ 
d}}}$\\
\\
\footnotesize{$^{\hbox{a }}$ School of Mathematics, University of Leeds, Leeds LS2 9JT, U.K.} \\ \footnotesize{Email: samuel.j.elliott@gmail.com}\\
\footnotesize{$^{\hbox{b }}$ INRIA Sophia Antipolis, 2004 route des Lucioles, BP 93, 06902 Sophia Antipolis, France.} \\ \footnotesize{Email: juliette.leblond@sophia.inria.fr}\\
\footnotesize{$^{\hbox{c }}$ Laboratoire Paul Painlev\'e, Cit\'e Scientifique Lille 1, 59655 Villeneuve-d'Ascq Cedex, France.} \\ \footnotesize{Email: elodie.pozzi@math.univ-lille1.fr}\\
\footnotesize{$^{\hbox{d }}$ Institut Fourier, 100 rue des maths, BP 74, 38402 Saint-Martin-d'H\`eres Cedex, France.} \\ \footnotesize{Email: emmanuel.russ@ujf-grenoble.fr}}
\maketitle


\begin{abstract}
Let $\Omega_1,\Omega_2\subset \CC$ be bounded domains. Let $\phi:\Omega_1\rightarrow \Omega_2$ holomorphic in $\Omega_1$ and belonging to $W^{1,\infty}_{\Omega_2}(\Omega_1)$. We study the composition operators $f\mapsto f\circ\phi$ on generalized Hardy spaces on $\Omega_2$, recently considered  in \cite{bfl, BLRR}. In particular, we provide necessary and/or sufficient conditions on $\phi$, depending on the geometry of the domains, ensuring that these operators are bounded, invertible, isometric or compact. Some of our results are new even for Hardy spaces of analytic functions.
\vskip 0.2cm
\noindent{\it Keywords:} Generalized Hardy spaces, composition operators.
\vskip 0.2cm
\noindent{\it MSC numbers:} Primary 47B33, secondary 30H10.
\end{abstract}

\tableofcontents

\section{Introduction}

The present work aims at generalizing properties of composition operators on Hardy spaces of domains of the complex plane to the framework of generalized Hardy spaces. 
Generalized analytic functions, among which pseudo-holomorphic functions, were  considered a long time ago,
see \cite{bn, vekua}, and more recently in \cite{krav}, in particular because of their links with classical partial differential equations (PDEs) in mathematical physics, like the conductivity or Schr\" odinger equations (see \cite[Lem. 2.1]{ap},  \cite{alessandrini-rondi}). 
By generalized analytic functions, we mean solutions (as distributions) to the following $\overline{\partial}$-type equations (real linear conjugate Beltrami and Schr\" odinger type elliptic PDEs): 
\[
\overline{\partial}f=\nu\overline{\partial f}\mbox{ or }\overline{\partial}w=\alpha\overline{w} \, ,
\]
without loss of generality (\cite{bn}).
For  specific classes of dilation coefficients $\nu$, $\alpha$, these two PDEs are equivalent to each other (as follows from a trick going back to Bers an Nirenberg, see  \cite{bn}). 
%
%
They are also related to the complex linear Beltrami equation, with the implicit dilation coefficient $\nu \overline{\partial f}/\partial f$, and to quasi--conformal applications \cite{ahlfors}.
Properties of associated (normed) Hardy classes $H^p_{\nu}$ and $G^p_{\alpha}$ have been established in  \cite{bfl,BLRR,EfendievRuss} for $1<p<\infty$ (these classes seem to have  been introduced in \cite{musaev} for simply connected domains). They
share many properties of the classical Hardy spaces of analytic (holomorphic) functions (for $\nu = \alpha = 0$).
The proofs of these properties rely on a  factorization result from \cite{bn} for generalized analytic functions which involve holomorphic functions. This factorization result was extended  in \cite{bbc,bfl,BLRR} to $G^p_{\alpha}$ functions, through classical Hardy spaces $H^p$, see Proposition \ref{hpgalpha}.

\noindent
Note that important applications of these classes come from Dirichlet--Neumann boundary value problems  and Cauchy type transmission issues for the elliptic conductivity PDE $\nabla \cdot \left( \sigma \, \nabla u \right) = 0$ with conductivity $\sigma = \left( 1-\nu\right) \, \left( 1+\nu\right)^{-1}$ in domains of $\RR^2 \simeq \CC$, see \cite{alessandrini-rondi, BLRR}. Indeed, on simply-connected domains, solutions $u$ coincide with real--parts of solutions $f$ to $\overline{\partial}f=\nu\overline{\partial f}$. In particular, this links Calder\'on's inverse
conductivity problem to similar issues for the real linear conjugate Beltrami equation, as in \cite{ap}.
Further, these new Hardy classes furnish a suitable framework in order to state and solve families of best constrained approximation issues (bounded extremal problems) \cite{EfendievRuss, flps}, from partial boundary values (given by Dirichlet--Neumann boundary conditions, through generalized harmonic conjugation or Hilbert transform). 

\noindent
In the Hilbertian setting $p=2$, constructive aspects are available for  particular conductivity coefficients $\nu$, for which bases of $H^2_\nu$ may be explicitly constructed, in the disk or the annulus, see \cite{fischer, flps}. In the annular setting, and in toroidal coordinates, this allows to tackle a free boundary problem related to  plasma confinment in tokamaks (\cite{fischer}). Namely, in toroidal plane sections, the boundary of the plasma is a level curve of the magnetic potential solution to a conductivity PDE. It has to be recovered from available magnetic data on the chamber (Dirichlet--Neumann data on the outer boundary of an annular domain). Bounded extremal problems provide a way to regularize and solve  this geometric inverse problem.
Observe that the boundary of the annular domain contained between the chamber and the plasma is not made of concentric 
circles. When it comes to realistic geometries, 
properties of composition operators on generalized Hardy classes may provide a selection of conformal maps from the disk or circular domains.\\


\noindent The present work is a study of composition operators on these Hardy classes. Let $\Omega\subset \CC$ be a domain. Hardy spaces $H_\nu^p(\Omega)$ of solutions to the conjugate Beltrami equation $\overline{\partial} f=\nu \overline{\partial f}$ {\it a.e} on $\Omega$ 
are first considered when $\Omega$ is the unit disc $\DD$ or in the annulus $\AAA=\{z\in\CC: r_0< \vert z\vert<1\}$. A way to define those spaces in bounded Dini-smooth domains (see below) is to use the conformal invariance property (see \cite{bfl}); more precisely, if $\Omega_1$ and $\Omega_2$ are two bounded Dini-smooth domains and $\phi$ a conformal map from $\Omega_1$ onto $\Omega_2$, then $f$ is in $H_\nu^p(\Omega_2)$, with $\nu\in W_{\RR}^{1,\infty}(\Omega_2)$ if and only if $f\circ\phi$ is in $H_{\nu\circ\phi}^p(\Omega_1)$ and $\nu\circ\phi\in W_{\RR}^{1,\infty}(\Omega_1)$. In terms of operator, if $\phi:\Omega_1\to\Omega_2 $ is an analytic conformal map, the composition operator $C_\phi:f\longmapsto f\circ\phi$ maps  $H_\nu^p(\Omega_2)$ onto $H_{\nu\circ\phi}^p(\Omega_1)$. Similar results hold in $G^p_{\alpha}$ Hardy spaces of solutions to $\overline{\partial}w=\alpha\overline{w}$. 

\noindent Suppose now that the composition map $\phi:\Omega_1\to\Omega_2$ is a function in $W^{1,\infty}(\Omega_1,\Omega_2)$ and analytic in $\Omega_1$, what can we say about $C_\phi(f)=f\circ\phi$ when $f\in H_\nu^p(\Omega_2)$ in terms of operator properties? This operator has been widely studied when $\Omega_1=\Omega_2=\DD$ and in the case of analytic (holomorphic) Hardy spaces $H^p(\DD)$ ({\it i.e.} $\nu\equiv 0$) giving characterizations of composition operators that are invertible in \cite{nord}, isometric  in \cite{Forelli},  similar to isometries in \cite{bayart}, and compact in \cite{shapiro1,shapsmith}, for example. Fewer  results are known concerning composition operators on $H^p$ spaces of an annulus. However, one can find in \cite{Boyd76} a sufficient condition on $\phi$ to have the boundedness of $C_\phi$ and a characterization of Hilbert-Schmidt composition operators. The study of composition operators has been generalized to many other spaces of analytic functions, such as Dirichlet spaces (\cite{llqr} and the references therein) or Bergman spaces (\cite{shapiro2}).\\ 

\noindent In this paper, we study some properties (boundedness, invertibility, isometry, compactness) for the composition operator defined on the Hardy space $H_\nu^p(\Omega)$ and $G^p_{\alpha}(\Omega)$ where $\Omega$ will be a bounded Dini-smooth domain (most of the time, $\Omega$ will be the unit disc $\DD$ or the annulus $\AAA$). 

In Section \ref{sec2}, we provide definitions of generalized Hardy classes together with some properties. 
Section \ref{sec:bdd} is devoted to boundedness results for composition operators on generalized Hardy classes for bounded Dini-smooth domains, while Section \ref{sec3} is related to their invertibility. Isometric composition operators on generalized Hardy classes of the disk and the annulus are studied in  Section \ref{sec:isom}, that appear to be new in $H^p(\AAA)$ as well. In Section \ref{sec:compact}, compactness properties for composition operators are investigated. A conclusion is written in Section \ref{sec:conclu}.
We will refer to specific results about analytic Hardy spaces $H^p$ thanks to factorization theorems (Appendix \ref{factor}), while properties of isometric composition operators on $H^p(\AAA)$ are established in Appendix \ref{app}.

\section{Definitions and notations}\label{sec2}

\subsection{Some notations}

\noindent In this paper, we will denote by $\Omega$ a connected open subset of the complex plane $\CC$ (also called a domain of $\CC$), by $\partial\Omega$ its boundary, by $\DD$ the unit disc and by $\TT=\partial\DD$ the unit circle. For $0<r_0<1$, let $\AAA$ be the annulus $\{z\in\CC :r_0<|z|<1\}=\DD\cap(\CC\backslash r_0\overline{\DD})$, the boundary of which is $\partial\AAA=\TT\cup\TT_{r_0}$, where $\TT_{r_0}$ is the circle of radius $r_0$. More generally, we will consider a circular domain $\GG$ defined as follows

\begin{equation}\label{circ_dom}
\GG=\DD\backslash\bigcup_{j=0}^{N-1} (a_j+r_j\overline{\DD}),
\end{equation}

\noindent where $N\geq 2$, $a_j\in\DD$, $0<r_j<1$, $0\leq j\leq N-1$. Its boundary is 

$$\partial\GG=\TT\cup\bigcup_{j=0}^{N-1} (a_j+\TT_{r_j})$$

\noindent where the circles $a_j+\TT_{r_j}$ for $0\leq j\leq N-1$ have a negative orientation whereas $\TT$ has the positive orientation. Note that for $N=2$ and $a_0=0$, $\GG$ is the annulus $\AAA$. \par
\noindent  A domain $\Omega$ of $\overline{\CC}$ is Dini-smooth if and only if its boundary $\partial\Omega$ is a finite union of Jordan curves with non-singular Dini-smooth parametrization. We recall that a function $f$ is said to be Dini-smooth if its derivative is Dini-continuous, \textit{i.e.} its modulus of continuity $\omega_f$ is such that

$$\int_0^{\varepsilon}\frac{\omega_f(t)}{t} dt<\infty,\,\,\hbox{for some}\,\,\varepsilon>0. $$
Recall that, if $\Omega$ is a bounded Dini-smooth domain, there exists a circular domain $\GG$ and a conformal map $\phi$ between $\GG$ and $\Omega$ which extends continuously to a homeomorphism between $\overline{\GG}$ and $\overline{\Omega}$, while the derivatives of $\phi$ also extend continuously to $\overline{\GG}$ (\cite{bfl}, Lemma A.1).
\noindent If $E,F$ are two Banach spaces, ${\cal{L}}(E,F)$ denotes the space of bounded linear maps from $E$ to $F$, and $T\in {\mathcal L}(E,F)$ is an isometry if and only if, for all $x\in E$, $\left\Vert Tx\right\Vert_F=\left\Vert x\right\Vert_E$.\\

\noindent If $A(f)$ and $B(f)$ are quantities depending on a function $f$ ranging in a set $E$, we will write $A(f)\lesssim B(f)$ when there is a positive constant $C$ such that $A(f)\leq C B(f)$ for all $f\in E$. We will say that $A(f)\sim B(f)$ if there is $C>0$ such that $C^{-1} B(f)\leq A(f)\leq C B(f)$ for all $f\in E$.

\subsection{Lebesgue and Sobolev spaces}

\noindent The Lebesgue measure on the complex plane will be denoted by $m$ and for a complex number $z=x+i y$

$$dm(z)=\frac{i}{2}\, dz\wedge d\overline{z}=dxdy. $$  

\noindent For $1\leq p<\infty$, $L^p(\Omega)$ designates the classical Lebesgue space of functions defined on $\Omega$ equipped with the norm

$$\Vert f\Vert_{L^p(\Omega)}:=\left(\int_{\Omega}\vert f(z)\vert^p dm(z)\right)^{1/p},\,\,f\in L^p(\Omega),$$
while $L^{\infty}(\Omega)$ stands for the space of essentially bounded measurable functions on $\Omega$ equipped with the norm
$$
\Vert f\Vert_{L^{\infty}(\Omega)}:=\mbox{ess sup}_{z\in \Omega} \left\vert f(z)\right\vert.
$$
\noindent We denote by ${\cal{D}}(\Omega)$ the space of smooth functions with compact support in $\Omega$.
Let ${\cal{D}}'(\Omega)$ be its dual space which is the space of distributions on $\Omega$.\\

\noindent For $1\leq p\leq \infty$, we recall that the Sobolev space $W^{1,p}(\Omega)$ is the space of all complex valued functions $f\in L^p(\Omega)$ with distributional derivatives in $L^p(\Omega)$. The space $W^{1,p}(\Omega)$ is equipped with the norm

$$\Vert f \Vert_{W^{1,p}(\Omega)}=\Vert f\Vert_{L^p(\Omega)}+\Vert\partial f\Vert_{L^p(\Omega)}+\Vert\overline{\partial} f\Vert_{L^p(\Omega)},$$ 

\noindent where the operators $\partial $ and $\overline{\partial}$ are defined, in the sense of distributions: for all $\phi\in{\cal{D}}(\Omega)$,

$$\langle \partial f,\phi\rangle=-\int f\partial\phi,\quad\langle \overline{\partial} f,\phi\rangle=-\int f\overline{\partial}\phi, $$

\noindent where

$$\partial=\frac{1}{2}(\partial_x-i\partial_y)\,\,\hbox{and}\,\,\overline{\partial}=\frac{1}{2}(\partial_x+i\partial_y). $$
Note that, when $\Omega$ is $C^1$ (in particular, when $\Omega$ is Dini-smooth), $W^{1,\infty}(\Omega)$ coincindes with the space of Lipschitz functions on $\Omega$ (\cite[Thm 4, Sec 5.8]{evans}). We will write $L_{\Omega_2}^p(\Omega)$ and $W_{\Omega_2}^{1,p}(\Omega)$ to specify that the functions have values in $\Omega_2\subset\CC$.  

\subsection{Hardy spaces}
For a detailed study of classical Hardy spaces of analytic functions, see \cite{duren,Garnett}. Let us briefly recall here some basic facts.\par 
\noindent For $1\leq p< \infty$, the Hardy space of the unit disc $H^p(\DD)$ is the collection of all analytic functions $f:\DD\longrightarrow\CC$ such that 

\begin{equation}\label{normD}
\Vert f\Vert_{H^p(\DD)}:=\left(\sup_{0<r<1}\frac{1}{2\pi}\int_{0}^{2\pi}\vert f(r e^{it})\vert^p dt\right)^{1/p}<\infty.
\end{equation}

\noindent For $p=\infty$, the Hardy space $H^{\infty}(\DD)$ is the Banach space of analytic functions which are bounded on $\DD$ equipped with the norm

$$\Vert f\Vert_{H^{\infty}(\DD)}=\sup_{z\in\DD}\vert f(z)\vert. $$

\noindent For $1\leq p\leq \infty$, any function $f\in H^p(\DD)$ has a non-tangential limit \textit{a.\,e.} on $\TT$ which we call the trace of $f$ and is denoted by $\mbox{tr}\,f$. For all $f\in H^p(\DD)$, we have that $\mbox{tr}\,f\in H^p(\TT)$ where $H^p(\TT)$ is a strict subspace of $L^p(\TT)$, namely

$$H^p(\TT)=\left\{h\in L^p(\TT),\hat{h}(n)=\frac{1}{2\pi}\int_{0}^{2\pi} h(e^{it})e^{-int} dt=0,\,\,n<0\right\}. $$  

\noindent More precisely, $H^p(\DD)$ is isomorphic to $H^p(\TT)$ and $\Vert f\Vert_{H^p(\DD)}=\Vert\mbox{tr}\,f\Vert_{L^p(\TT)}$, which allows us to identify the two spaces $H^p(\DD)$ and $H^p(\TT)$.\\

\noindent Likewise, in \cite{Sarason65}, the Hardy space $H^p(\AAA)$ of an annulus $\AAA$ is the space of analytic functions $f$ on $\AAA$ such that

\begin{equation}\label{normA}
\Vert f\Vert_{H^p(\AAA)}=\left(\sup_{r_0<r<1}\frac{1}{2\pi}\int_{0}^{2\pi}\vert f(r e^{it})\vert^p dt\right)^{1/p}<\infty,
\end{equation}

\noindent for $1\leq p<\infty$. It can also be viewed as the topological direct sum

\begin{equation}\label{annulus:direct}
H^p(\AAA)=H^p(\DD)_{|_\AAA}\oplus H^p(\overline{\CC}\backslash r_0\overline{\DD})_{|_{\AAA}}, 
\end{equation}

\noindent where $\overline{\CC}=\CC\cup\{\infty\}$ and $H^p(\overline{\CC}\backslash r_0\overline{\DD})$ is isometrically isomorphic to $H^p(r_0\DD)$ \textit{via} the transformation

$$f\in H^p(\DD)\longmapsto \widetilde{f}\in H^p(\overline{\CC}\backslash r_0\overline{\DD}),$$

\noindent where $\widetilde{f}(z)=f\left(\frac{r_0}{z}\right)$ for all $z\in\overline{\CC}\backslash r_0\overline{\DD}$. It follows from \eqref{annulus:direct} that any function $f\in H^p(\AAA)$ has a non-tangential limit \textit{a.\,e.} on $\partial\AAA$ also denoted by $\mbox{tr}\,f$, such that 

$$\mbox{tr}\,f\in H^p(\partial\AAA)=\left\{h\in L^p(\partial\AAA),\widehat{h_{|_{\TT}}}(n)=r_0^n\,\widehat{h_{|_{r_0\TT}}}(n),\,\,n\in\ZZ\right\}.$$

\noindent where for $n\in\ZZ$

$$\widehat{h_{|_{\TT}}}(n)=\frac{1}{2\pi}\int_{0}^{2\pi} h(e^{it})e^{-int} dt $$

\noindent and 

$$ \widehat{h_{|_{r_0\TT}}}(n)=\frac{1}{2\pi}\int_{0}^{2\pi} h(r_0\,e^{it})e^{-int}.$$

\noindent are respectively the $n$-th Fourier coefficients of $h_{|_{\TT}}$ and $h_{|_{r_0\TT}}$. Again, the space $H^p(\AAA)$ can be identified to $H^p(\partial\AAA)$ \textit{via} the isomorphic isomorphism $f\in H^p(\AAA)\longmapsto \mbox{tr}\,f\in H^p(\partial\AAA)$ and thus $\Vert f\Vert_{H^p(\AAA)}=\Vert\mbox{tr}\,f\Vert_{L^p(\partial\AAA)}$.\\

\noindent The definition of Hardy spaces has been extended in \cite{rudin} to any complex domain $\Omega$ using harmonic majorants. More precisely, for $1\leq p<\infty$ and $z_0\in\Omega$, $H^p(\Omega)$ is the space of analytic functions $f$ on $\Omega$ such that there exists a harmonic function $u:\Omega\longrightarrow [0,\infty)$ such that for $z\in\Omega$

$$\vert f(z)\vert^p\leq u(z).$$ 

\noindent The space is equipped with the norm

$$\inf\left\{u(z_0)^{1/p},\,\,\vert f\vert^p\leq u\,\,\hbox{for}\,\,u\,\,\hbox{harmonic function in}\,\,\Omega\right\}$$ 

\begin{remark}\label{rmk1}
\begin{itemize}
\item[$1.$]
It follows from the Harnack inequality (\cite[3.6, Ch.3]{abr}) that different choices of $z_0$  give rise to equivalent norms in $H^p(\Omega)$.
\item[$2.$]
\noindent If $\Omega=\DD$ or $\AAA$, the two previously defined norms on $H^p(\Omega)$ are equivalent.
\end{itemize}
\end{remark}

%

\subsection{Generalized Hardy spaces}

\subsubsection{Definitions} \label{defhardygen}

\noindent Let $1<p< \infty$ and $\nu\in W^{1,\infty}_{\RR}(\DD)$ such that $\Vert \nu\Vert_{L^{\infty}(\DD)}\leq\kappa$ with $\kappa\in (0,1)$. The generalized Hardy space of the unit disc $H_{\nu}^p(\DD)$ was first defined in \cite{musaev} and then in \cite{BLRR} as the collection of all measurable functions $f:\DD\longrightarrow\CC$ such that $\overline{\partial}f=\nu\,\overline{\partial f}$ in the sense of distributions in $\DD$ and

\begin{equation}\label{norm_nuD}
\Vert f\Vert_{H_{\nu}^p(\DD)}:=\displaystyle{\left(\essup_{0<r<1}\frac{1}{2\pi}\int_{0}^{2\pi}\vert f(re^{it})\vert^p dt\right)^{1/p}<\infty}.
\end{equation}
\noindent The definition was extended, in \cite{EfendievRuss}, to the annulus $\AAA$: for $\nu\in W_{\RR}^{1,r}(\AAA)$, $r\in (2,\infty)$, $H_{\nu}^p(\AAA)$ is the space of functions $f:\AAA\to\CC$ such that $\overline{\partial}f=\nu\,\overline{\partial f}$ in the sense of distribution in $\AAA$ and satisfying

\begin{equation}\label{norm_nuA}
\displaystyle
\Vert f\Vert_{H_{\nu}^p(\AAA)}:=\left(\essup_{r_0<r<1}\frac{1}{2\pi}\int_{0}^{2\pi}\vert f(re^{it})\vert^p dt\right)^{1/p}<\infty.
\end{equation}

\noindent Now, let $\Omega\subset \CC$ be a domain and $\nu$ such that

\begin{equation}
\label{kappa}
\nu\in W^{1,\infty}_{\RR}(\Omega) \, , \ \Vert \nu\Vert_{L^{\infty}(\Omega)}\leq\kappa \, , \mbox{ with } \kappa\in(0,1) \, .  
\end{equation}

\noindent The definition of $H^p_{\nu}(\Omega)$ was further extended to the case where $\Omega$ is a Dini-smooth domain of $\overline{\CC}$ (see \cite{bfl}). In this case, the norm is defined by
\begin{equation} \label{Smirnovnorm}
\left\Vert g\right\Vert_{H^p_{\nu}(\Omega)}:=\sup_{n\in\NN}\Vert g\Vert_{L^p(\partial\Delta_n)},
\end{equation}
where  $(\Delta_n)_n$ is a fixed sequence of domains such that $\overline{\Delta_n}\subset\Omega$ and $\partial\Delta_n$ is a finite union of rectifiable Jordan curves of uniformly bounded length, such that each compact subset of $\Omega$
 is eventually contained in $\Delta_n$ for $n$ large enough. We refer to \cite{bfl} for the existence of such sequence. \\

\noindent In parallel with Hardy spaces $H_{\nu}^p(\Omega)$ (with $\Omega$ equal to $\DD$, $\AAA$ or more generally to a Dini-smooth domain), Hardy spaces $G_\alpha^p(\Omega)$ were defined in \cite{bfl,BLRR,EfendievRuss,musaev} for $\alpha\in L^{\infty}(\Omega)$ as the collection of measurable functions $w:\Omega\to\CC$ such that $\overline{\partial}w=\alpha\,\overline{w}$ in ${\mathcal D}^{\prime}(\Omega)$ and 

\begin{equation}\label{norm_alpha_A}
\Vert w\Vert_{G_{\alpha}^p(\Omega)}=\left({\essup}_{\rho<r<1}\frac{1}{2\pi}\int_{0}^{2\pi}\vert w(re^{it})\vert^p dt\right)^{1/p}<\infty,
\end{equation}
 
\noindent with $\rho=0$ if $\Omega=\DD$ and $\rho=r_0$ if $\Omega=\AAA$. If $\Omega$ is a Dini-smooth domain, the essential supremum is taken over all the $L^p(\partial\Delta_n)$ norm of $w$ for $n\in\NN$.

\begin{remark}
The generalized Hardy spaces $H_{\nu}^p(\Omega)$ and $G_{\alpha}^p(\Omega)$ are real Banach spaces (note that when $\nu=0$ or $\alpha=0$ respectively, they are complex Banach spaces). 
\end{remark}
\noindent Recall that if $\Omega$ is a bounded Dini-smooth domain, a function $g$ lying in generalized Hardy spaces $H_\nu^p(\Omega)$ or $G_{\alpha}^p(\Omega)$ has a non-tangential limit \textit{a.e.} on $\partial\Omega$ which is called the trace of $g$ is denoted by $\mbox{tr}\,g\in L^p(\partial\Omega)$ and 

\begin{equation}\label{normtrace}
\Vert g\Vert_{H_\nu^p(\Omega)}\sim \Vert\mbox{tr}\,g\Vert_{L^p(\partial\Omega)}, 
\end{equation}

\noindent (see \cite{bfl,BLRR,EfendievRuss}). We will denote by $\mbox{tr}\,(H^p_{\nu}(\Omega))$ the space of traces of $H_\nu^p(\Omega)$-functions; it is a strict subspace of $L^p(\partial\Omega)$. Note also that $g\mapsto \Vert\mbox{tr}\,g\Vert_{L^p(\partial\Omega)}$ is a norm on $H^p_{\nu}(\Omega)$, equivalent to the one given by \eqref{Smirnovnorm}. However, contrary to the case of Hardy spaces of analytic functions of the disk, $\left\Vert \cdot \right\Vert_{H^p_{\nu}(\DD)}$ and $\left\Vert \mbox{ tr } \cdot \right\Vert_{L^p(\TT)}$ are not equal in general (see \eqref{normtrace}).\par
\noindent Finally, functions in $H^p_{\nu}(\Omega)$ and $G^p_{\alpha}(\Omega)$ are continuous in $\Omega$:
\begin{lemma} \label{contomega}
Let $\Omega\subset \CC$ be a bounded Dini-smooth domain, $\nu\in W^{1,\infty}_{\RR}(\Omega)$ meeting \eqref{kappa} and $\alpha\in L^{\infty}(\Omega)$. Then, all functions in $G^p_{\alpha}(\Omega)$ and $H^p_{\nu}(\Omega)$ are continuous in $\Omega$.
\end{lemma}
\beginpf
Indeed, let $\omega\in G^p_{\alpha}(\Omega)$. By \cite[Prop. 3.2]{bfl}, $\omega=e^sF$ with $s\in C(\overline{\Omega})$ (since $s\in W^{1,r}(\Omega)$ for some $r>2$) and $F\in H^p(\Omega)$. Thus, $\omega$ is continuous in $\Omega$. If $f\in H^p_{\nu}(\Omega)$ and $\omega=\mathcal{J}^{-1}(g)$, then $\omega\in G^p_{\alpha}(\Omega)$ is continuous in $\Omega$ and since $\nu$ is continuous and \eqref{kappa} holds, $f\in C(\Omega)$.

\subsubsection{An equivalent norm}

\noindent Throughout the present section, unless explicitly stated, let $\Omega$ be an arbitrary bounded domain of $\CC$. For $1<p<\infty$, we define generalized Hardy spaces on $\Omega$, inspired by the definitions of Hardy spaces of analytic functions given in \cite{rudin}. Let $\nu$ meet \eqref{kappa}.
\begin{definition} \label{defhpnu}
Define $E^p_{\nu}(\Omega)$ as the space of measurable functions $f:\Omega\rightarrow \CC$ solving 

\begin{equation} \label{CB}
\overline{\partial}f=\nu\,\overline{\partial f}\mbox{ in }{\mathcal D}^{\prime}(\Omega),
\end{equation}
and for which there exists a harmonic function $u:\Omega\rightarrow [0,+\infty)$ such that
\begin{equation} \label{majorharm}
\left\vert f(z)\right\vert^p\leq u(z)
\end{equation}
for almost every $z\in \Omega$. Fix a point $z_0\in \Omega$ and define

\begin{equation}\label{normharmH}
\left\Vert f\right\Vert_{E^p_{\nu}(\Omega)}:=\inf u^{1/p}(z_0),
\end{equation}

the infimum being taken over all harmonic functions $u:\Omega\rightarrow [0,+\infty)$ such that \eqref{majorharm} holds.
\end{definition}

\noindent Let $\alpha\in L^{\infty}(\Omega)$. Let us similarly define $F^p_{\alpha}(\Omega)$: 
\begin{definition}\label{defgalpha}
Define $F^p_{\alpha}(\Omega)$ as the space of measurable functions $w:\Omega\rightarrow \CC$ solving
\begin{equation} \label{eqalpha}
\overline{\partial}w=\alpha\,\overline{w}\mbox{ in }{\mathcal D}^{\prime}(\Omega),
\end{equation}
and for which there exists a harmonic function $u:\Omega\rightarrow [0,+\infty)$ such that \eqref{majorharm} holds. Define $\left\Vert w\right\Vert_{F^p_{\alpha}(\Omega)}$ by

\begin{equation}\label{normharmG}
\left\Vert w\right\Vert_{F^p_{\alpha}(\Omega)}:=\inf u^{1/p}(z_0),
\end{equation}

where the infimum is computed as in Definition \ref{defhpnu}.
\end{definition}

\noindent Observe that in the above definitions, different values of $z_0$ give rise to equivalent norms as in Remark \ref{rmk1}. We first check:
\begin{proposition}
\begin{itemize}
\item[$1.$]
The map $f\mapsto \left\Vert f\right\Vert_{E^p_{\nu}(\Omega)}$ is a norm on $E^p_{\nu}(\Omega)$.
\item[$2.$]
The analogous conclusion holds for $F^p_{\alpha}(\Omega)$.
\end{itemize}
\end{proposition}
\beginpf
It is plain to see that $\left\Vert \cdot\right\Vert_{E^p_{\nu}(\Omega)}$ is positively homogeneous of degree $1$ and subadditive. Assume now that $\left\Vert f\right\Vert_{H^p_{\nu}(\Omega)}=0$. That $f=0$ follows at once from the fact that, if $(u_j)_{j\geq 1}$ is a sequence of nonnegative harmonic functions on $\Omega$ such that $u_j(z_0)\rightarrow 0$, $j \to \infty$, for $z_0 \in \Omega$ from Definition \ref{defhpnu}, then $u_j(z)\rightarrow 0$, $j \to \infty$,  for all $z\in \Omega$. To check this fact, define
$$
A:=\left\{z\in \Omega;\ u_j(z)\rightarrow 0\right\}.
$$
The Harnack inequality (\cite[3.6, Ch.3]{abr}) shows at once that $A$ is open in $\Omega$. If $B=\Omega\setminus A$, then the Harnack inequality also shows that $B$ is open. Because $z_0 \in A\neq \emptyset$ and $\Omega$ is connected, then $A=\Omega$, which proves point 1 and, similarly, point 2. 
\endpf

%

\noindent Let $\nu$ satisfying assumption (\ref{kappa}) and $\alpha\in L^{\infty}(\Omega)$ associated with $\nu$ in the sense that
\begin{equation}\label{nu-lambda}
\alpha=\frac{- \overline{\partial}\nu}{{1-\nu^2}},
\end{equation}
\noindent The link between $E^p_{\nu}(\Omega)$ and $F^p_{\alpha}(\Omega)$ is as follows (see \cite{bfl,BLRR} in the case of Dini-smooth domains):
\begin{proposition} \label{hpgalpha}
A function $f:\Omega\rightarrow \CC$ belongs to $E^p_{\nu}(\Omega)$ if and only if
\begin{equation} \label{wf}
w ={\mathcal J}(f):=\frac{f-\nu\overline{f}}{\sqrt{1-\nu^2}}
\end{equation}
belongs to $F^p_{\alpha}(\Omega)$. One has $\left\Vert f\right\Vert_{E^p_{\nu}(\Omega)}\sim \left\Vert w\right\Vert_{F^p_{\alpha}(\Omega)}$.
\end{proposition}

\beginpf That $f$ solves \eqref{CB} if and only if $w$ solves \eqref{eqalpha} was checked in \cite{bfl,BLRR}. That $\left\vert f\right\vert^p$ has a harmonic majorant if and only if the same holds for $\left\vert w\right\vert^p$ and $\left\Vert f\right\Vert_{E^p_{\nu}(\Omega)}\sim \left\Vert w\right\Vert_{F^p_{\alpha}(\Omega)}$ are straightforward consequences of \eqref{wf} and assumption \eqref{kappa}.
\endpf\par

\begin{remark}\label{equiv:def}
As \cite[Thm 3.5, $(ii)$]{bfl} shows, when $\Omega$ is a Dini-smooth domain, $\nu$ meets \eqref{kappa} and $\alpha\in L^{\infty}(\Omega)$, $H^p_{\nu}(\Omega)=E^p_{\nu}(\Omega)$ and $G^p_{\alpha}(\Omega)=F^p_{\alpha}(\Omega)$, with equivalent norms. In this case, if $\Omega$ is Dini-smooth, then, for $w\in G_\alpha^p(\Omega)$ we have that

$$\Vert w\Vert_{G_\alpha^p(\Omega)}\sim\inf u^{1/p}(z_0),$$

\noindent where the infimum is taken as in Definition \ref{defgalpha}. The same stands for $f\in H_\nu^p(\Omega)$.  
\end{remark} 

\noindent Proposition \ref{hpgalpha} immediately yields:
\begin{lemma}\label{op-nu-alpha}
Let $\nu,\tilde{\nu}$ satisfying \eqref{kappa} and $\alpha,\tilde{\alpha}$ associated with $\nu$ (resp. $\tilde{\nu}$) as in equation \eqref{nu-lambda}. Then, $T\in{\cal{L}}(H_\nu^p(\Omega),H_{\widetilde{\nu}}^p(\Omega))$ if and only if $\widetilde{T}\in{\cal{L}}(G_\alpha^p(\Omega),G_{\widetilde{\alpha}}^p(\Omega))$ where $\widetilde{\Cj} T=\widetilde{T}\cal{J}$, 
and $\widetilde{\Cj}$ is the 
$\RR$-linear isomorphism from $H_{\widetilde{\nu}}^p(\DD)$ onto $G_{\widetilde{\alpha}}^p(\Omega)$ 
defined by (\ref{wf}) with $\nu$ replaced by $\widetilde{\nu}$.
\end{lemma}

\section{Boundedness of composition operators on generalized Hardy spaces}
\label{sec:bdd}

\noindent Let $\Omega_1,\Omega_2$ be two bounded Dini-smooth domains in $\CC$, $\nu $ defined on  $\Omega_2$ satisfying assumption (\ref{kappa}), 
and $\phi$ satisfying:
\begin{equation}
\label{phi}
\phi:\Omega_1\to\Omega_2 \mbox{ analytic with } \phi\in W^{1,\infty}_{\Omega_2}(\Omega_1). 
\end{equation}


We consider the composition operator $C_\phi$ defined on $H_\nu^p(\Omega_2)$ by $C_\phi(f)=f\circ\phi$.\par
\noindent Observe first that $\nu\circ\phi\in W^{1,\infty}_{\RR}(\Omega_1)$ since $\nu$ and $\phi$ are Lipschitz functions in $\Omega_2$ and $\Omega_1$ respectively and $\left\Vert \nu\circ\phi\right\Vert_{L^{\infty}(\Omega_1)}\leq \kappa$; hence $\nu\circ\phi$ satisfies (\ref{kappa}) on $\Omega_1$.
\begin{proposition} \label{continuouscphi}
The composition operator $C_\phi:H_\nu^p(\Omega_2)\to H_{\nu\circ\phi}^p(\Omega_1)$ is continuous.
\end{proposition}
\beginpf Let $f\in H^p_\nu(\Omega_2)$. Observe that $f\circ\phi$ is a Lebesgue measurable function on $\Omega_1$ and, since $\overline{\partial}\phi=0$ in $\Omega_1$,
\begin{eqnarray*}
\dbar (f\circ\phi)&=&[(\dbar f)\circ\phi]\dbar(\overline{\phi})=(\nu\circ\phi)[\overline{\partial f}\circ\phi]\dbar (\overline{\phi})\\
& &\\
                             &=&(\nu\circ\phi)\overline{(\partial f\circ\phi)\partial\phi}=(\nu\circ\phi)\overline{\partial (f\circ\phi)},
\end{eqnarray*}
\noindent (equalities are considered in the sense of distributions). Now, if $u$ is any harmonic majorant of $\vert f\vert^p$ in $\Omega_2$, then $u\circ\phi$ is a harmonic majorant of $\vert f\circ\phi\vert^p$ in $\Omega_1$, which proves that $C_\phi(f)\in H_{\nu\circ\phi}^p(\Omega_1)$. Moreover, by the Harnack inequality applied in $\Omega_2$,
$$
\left\Vert C_\phi(f)\right\Vert_{H^p_{\nu\circ\phi}(\Omega_1)}\leq u(\phi(z_0))^{1/p}\leq C \, u(z_0)^{1/p},
$$
for $z_0 \in \Omega_2$ as in Definition \ref{defhpnu}, and where the constant $C$ depends on $\Omega_2$, $z_0$ and $\phi(z_0)$ but not on $u$, so that, taking the infimum over all harmonic functions $u\geq \left\vert f\right\vert^p$ in $\Omega_2$, one concludes
$$
\left\Vert C_\phi(f)\right\Vert_{H^p_{\nu\circ\phi}(\Omega_1)}\lesssim\left\Vert f\right\Vert_{H^p_{\nu}(\Omega_2)}.
$$
\endpf
\begin{remark}
In the case where $\Omega_1=\Omega_2=\DD$, if $H^p_{\nu}(\DD)$ and $H^p_{\nu\circ\phi}(\DD)$ are equipped with the norms given by \eqref{normharmH}, the following upper bound for the operator norm of $C_\phi$ holds:
$$
\Vert C_\phi\Vert\leq\left(\frac{1+\vert\phi(0)\vert}{1-\vert\phi(0)\vert}\right)^{1/p}.
$$
Indeed, if $u$ is as before, one obtains
\begin{eqnarray*}
u\circ\phi(0)&=&\frac{1}{2\pi}\int_0^{2\pi}\frac{1-\vert\phi(0)\vert^2}{\vert e^{it}-\phi(0)\vert^2}u(e^{it}) dt\\
                     &\leq&\frac{1+\vert\phi(0)\vert}{1-\vert\phi(0)\vert}\frac{1}{2\pi}\int_0^{2\pi} u(e^{it}) dt=\frac{1+\vert\phi(0)\vert}{1-\vert\phi(0)\vert}u(0).
\end{eqnarray*}\par
\noindent In the doubly-connected case, assume that $\Omega=\AAA$. Let $z_0\in \AAA$ and $\psi$ be an analytic function from $\DD$ onto $\AAA$ such that $\psi(0)=z_0$. Arguing as in \cite{Boyd76}, we obtain an ``explicit'' upper bound for $\left\Vert C_{\phi}\right\Vert$. Indeed, let $u$ be as before. Using the harmonicity of $u\circ\psi$ in $\DD$, for all $s$ such that $\psi(s)=\phi(z_0)$, one has, for all $r\in (\left\vert s\right\vert,1)$,
\begin{eqnarray*}
u(\phi(z_0))=u(\psi(s))&=&\frac{1}{2\pi}\int_{0}^{2\pi}\Re\left(\frac{r\,e^{it}+s}{r\,e^{it}-s}\right)u\circ\psi(r\,e^{it}) dt\\
                                 &\leq&\frac{r+\vert s\vert}{r-\vert s\vert} u(\psi(0))=\frac{r+\vert s\vert}{r-\vert s\vert} u(z_0).
\end{eqnarray*}

\noindent Letting $r$ tend to $1$, we obtain 

$$u(\phi(z_0))\leq \inf_{s\in\psi^{-1}(\phi(z_0))} \frac{1+\vert s\vert}{1-\vert s\vert}.\,u(z_0),$$

\noindent which, with Definition \ref{defhpnu}, yields $\Vert C_\phi\Vert\leq\left(\displaystyle{\inf_{s\in\psi^{-1}(\phi(z_0))} \frac{1+\vert s\vert}{1-\vert s\vert}}\right)^{1/p}$. 
\end{remark}


\begin{remark}
Note that the conclusion of Proposition \ref{continuouscphi} and its proof remain valid when $\Omega_1$ and $\Omega_2$ are arbitrary connected open subsets of $\CC$.  
\end{remark}

\noindent In the sequel, when necessary, we will consider the composition operator defined on $G_\alpha^p$ spaces instead of $H^p_{\nu}$ spaces.  The next lemma shows that a composition operator defined on $H_\nu^p$ spaces is $\RR$-isomorphic to a composition operator on $G_\alpha^p$ spaces.

\begin{lemma} \label{boundedgpalpha}
Let $\nu$ (resp. $\phi$) satisfying (\ref{kappa}) (resp. (\ref{phi})). The composition operator $C_\phi$ mapping $H^p_{\nu}(\Omega_2)$ to $H^p_{\widetilde{\nu}}(\Omega_1)$ with $\widetilde{\nu}=\nu\circ\phi$ is then equivalent to the composition operator 
$\widetilde{C_\phi}$ mapping $G^p_{\alpha}(\Omega_2)$ to $G^p_{\widetilde{\alpha}}(\Omega_1)$, where $\widetilde{\alpha}$ is associated with $\widetilde{\nu}$ through (\ref{nu-lambda}). 
Moreover, 
\begin{equation}
\label{alpha}
\widetilde{\alpha}=(\alpha\circ\phi)\overline{\partial\phi} \, .
\end{equation}
In other words, for $\Cj$, $\widetilde{\Cj}$ defined as in Lemma \ref{op-nu-alpha},
we have the following commutative diagram:
$$\begin{CD}   
H_{\nu}^p(\Omega_2)@>{C_{\phi}}>{}>H_{\widetilde{\nu}}^p(\Omega_1)\\
@V{\Cj}V{}V@V{\widetilde{\Cj}}V{}V\\
G_{\alpha}^p(\Omega_2)@>{\widetilde{C}_{\phi}}>{}>G_{\widetilde{\alpha}}^p(\Omega_1)
\end{CD}$$
\end{lemma}
%
%
\beginpf The inverse of $\Cj$ is given by (see \cite{BLRR}):
\begin{equation}
\label{i-1}
\Cj^{-1} :w\in G_{\alpha}^p(\Omega_2)\longmapsto f=\frac{w+\nu\,\overline{w}}{\sqrt{1-\nu^2}}\in H_{\nu}^{p}(\Omega_2). 
\end{equation}
\noindent Note that 
$$\widetilde{\alpha}= \frac{- \dbar\widetilde{\nu}}{{1-\widetilde{\nu}^2}}=\frac{-\dbar(\nu\circ\phi)}{{1-\nu^2\circ\phi}}=\frac{-[(\dbar\nu)\circ\phi]\overline{\partial\phi}}{{1-\nu^2\circ\phi}}=(\alpha\circ\phi)\overline{\partial\phi},$$
\noindent  and $\widetilde{\Cj}$ is also an $\RR$-linear isomorphism from $H^p_{\widetilde{\nu}}(\Omega_1)$ onto $G^p_{\widetilde{\alpha}}(\Omega_1)$. Now, for any $f\in H^p_{\nu}(\Omega_2)$, we have that 
$$\widetilde{\Cj}(C_{\phi}(f))=\frac{f\circ\phi-(\nu\circ\phi)\overline{f\circ\phi}}{\sqrt{1-\nu^2\circ\phi}}=\left[\frac{f-\nu\overline{f}}{\sqrt{1-\nu^2}}\right]\circ\phi=\widetilde{C_\phi}(\Cj(f)).$$
%
%
%

\endpf

\section{Invertibility of the composition operator on $H_\nu^p(\Omega)$}\label{sec3}

%
%

\noindent In this section, we characterize invertible composition operators between $H^p_{\nu}$ spaces. \par
\noindent We will need an observation on the extension of a function $\nu$ meeting condition \eqref{kappa}. Before stating it, let us recall that, if $\Omega_1$ and $\Omega_2$ are open subsets of $\CC$, the notation $\Omega_1\subset\subset \Omega_2$ means that $\overline{\Omega_1}$ is a compact included in $\Omega_2$.
\begin{lemma} \label{extension}
Let $\Omega_1\subset\subset \Omega_2\subset \CC$ be bounded domains and $\nu$ be a Lipschitz function on $\Omega_1$ meeting condition \eqref{kappa}. There exists a Lipschitz function $\widetilde{\nu}$ on $\CC$ such that:
\begin{itemize}
\item[$1.$]
$\widetilde{\nu}(z)=\nu(z)$ for all $z\in \Omega_1$,
\item[$2.$]
the support of $\widetilde{\nu}$ is a compact included in $\Omega_2$,
\item[$3.$]
$\left\Vert \widetilde{\nu}\right\Vert_{L^{\infty}(\CC)}<1$.
\end{itemize}
\end{lemma}
\beginpf
Extend first $\nu$ to a compactly supported Lipschitz function on $\CC$, denoted by $\nu_1$. There exists an open set $\Omega_3$ such that $\Omega_1\subset\subset \Omega_3\subset\subset \Omega_2$ and $\left\Vert \nu_1\right\Vert_{L^{\infty}(\Omega_3)}<1$. Let $\chi\in {\mathcal D}(\CC)$ be such that $0\leq \chi(z)\leq 1$ for all $z\in \CC$, $\chi(z)=1$ for all $z\in \Omega_1$ and $\chi(z)=0$ for all $z\notin \Omega_3$. The function $\widetilde{\nu}:=\chi\nu_1$ satisfies all the requirements.
\endpf

Let $1<p<+\infty$, $\Omega\subset \CC$ be a bounded Dini-smooth domain and $\nu$ meet \eqref{kappa}.  For $z\in \Omega$, let ${\mathcal E}^{\nu}_z$, ${\mathcal F}^{\nu}_z$ be the real-valued evaluation maps at $z$ defined on $H^p_{\nu}(\Omega)$ and $G_\alpha^p(\Omega)$ 
by
$$
{\mathcal E}^{\nu}_z(f):=\mbox{ Re }f(z)\mbox{ and } {\mathcal F}^{\nu}_z(f):=\mbox{ Im }f(z) \mbox{ for all } f \in H^p_{\nu}(\Omega) \, , \ f \in G_\alpha^p(\Omega) \, . 
$$
\begin{proposition} \label{evalbounded}
For $z\in\Omega$, the evaluation maps ${\mathcal E}^{\nu}_z$ and ${\mathcal F}^{\nu}_z$ are continuous on $H_{\nu}^p(\Omega)$ and $G_\alpha^p(\Omega)$. 
\end{proposition}
\beginpf 
Let $f\in H^p_{\nu}(\Omega)$ and $z\in\Omega$. By definition \ref{defhpnu} of the norm in $H^p_{\nu}(\Omega)$, there exists a harmonic function $u$ in $\Omega$ such that $\left\vert f\right\vert^p\leq u$ in $\Omega$ with $u^{1/p}(z_0)\lesssim \left\Vert f\right\Vert_{H^p_{\nu}(\Omega)}$ for a fixed $z_0\in\Omega$. The Harnack inequality then yields
$$
\left\vert f(z)\right\vert\leq u^{1/p}(z)\lesssim u^{1/p}(z_0)\lesssim \Vert f\Vert_{H_{\nu}^p(\Omega)}$$

\noindent and thus we have

$$\vert\mbox{Re }f(z)\vert\lesssim\left\Vert f\right\Vert_{H^p_{\nu}(\Omega)} \mbox{ and } \vert\mbox{Im }f(z)\vert\lesssim\left\Vert f\right\Vert_{H^p_{\nu}(\Omega)},$$ 
which ends the proof. 

\endpf\par

%
%
%

\noindent For the characterization of invertible composition operators on $H^p_{\nu}$ spaces, we will need the fact that $H^p_{\nu}(\Omega)$ separates points in $\Omega$, when $\Omega$ is a Dini-smooth domain:

\begin{lemma} \label{separation}
Assume that $\Omega\subset \CC$ is a bounded Dini-smooth domain. Let $z_1\neq z_2\in \Omega$. Then, there exists $f\in H^p_{\nu}(\Omega)$ such that $f(z_1)\neq f(z_2)$.
\end{lemma}
\beginpf 
There exists $F\in H^p(\Omega)$ such that $F(z_1)=0$ and $F(z_2)\neq 0$ (take, for instance, $F(z)=z-z_1$). By Theorem \ref{factorgpalpha} in Appendix \ref{factor} below, there exists $s\in W^{1,r}(\Omega)$, for some $r\in (2,+\infty)$ such that $w=e^sF\in G^p_{\alpha}(\Omega)$. One has $w(z_1)=0$ and $w(z_2)\neq 0$. If $f:={\mathcal J}^{-1}(w)=\frac{w+\nu\overline{w}}{\sqrt{1-\nu^2}}$, $f\in H^p_{\nu}(\Omega)$ by Proposition \ref{hpgalpha}, $f(z_1)=0$ and $f(z_2)\neq 0$, since $\left\Vert \nu\right\Vert_{L^{\infty}(\Omega)}<1$. 
\endpf\par
We will also use in the sequel a regularity result for a solution of a Dirichlet problem for equation \eqref{CB}, where the boundary data is $C^{1}$ and only prescribed on one curve of $\partial\Omega$ :
\begin{lemma} \label{solsobolev}
Let $\Omega\subset \CC$ be a bounded $n$-connected Dini-smooth domain. Write $\partial\Omega=\cup_{j=0}^n \Gamma_j$, where the $\Gamma_j$ are pairwise disjoint Jordan curves. Fix $j\in \left\{0,...,n\right\}$. Let $\nu$ meet \eqref{kappa} and $\psi\in C^{1}_{\RR}(\Gamma_j)$. There exists $f\in H^{p}_{\nu}(\Omega)$ such that $\mbox{Re tr }f=\psi$ on $\Gamma_j$ and $\left\Vert f\right\Vert_{H^p_{\nu}(\Omega)}\lesssim\left\Vert \psi\right\Vert_{L^p(\Gamma_j)}$. Moreover, $f\in C(\overline{\Omega})$. 
\end{lemma}
\beginpf {\bf Step 1: } Let us first assume that $\Omega=\DD$. Since $\psi\in W^{1-1/q,q}_{\RR}(\TT)$ for some $q>\max(2,p)$, the result \cite[Thm 4.1.1]{BLRR} shows that there exists $f\in W^{1,q}(\DD)$ solving $\overline{\partial}f=\nu\overline{\partial f}$ in $\DD$ with $\mbox{Re tr }f=\psi$ on $\TT$. By \cite[Prop. 4.3.3]{BLRR}, $f\in H^q_{\nu}(\DD)\subset H^p_{\nu}(\DD)$, and since $q>2$, $f$ is continuous on $\overline{\DD}$. \par
\noindent  {\bf Step 2: }ÊAssume that $\Omega=\CC\setminus r_0\overline{\DD}$ for some $r_0\in (0,1)$. Let $\psi\in C^1_{\RR}(r_0\TT)$. For all $z\in \TT$, define $\widetilde{\psi}(z):=\psi\left(\frac{r_0}{\overline{z}}\right)$ and, for all $z\in \DD$, define $\widetilde{\nu}(z):=\nu\left(\frac{r_0}{\overline{z}}\right)$. Step $1$ yields a function $\widetilde{f}\in H^p_{\widetilde{\nu}}(\DD)$, continuous on $\overline{\DD}$, such that $\mbox{Re tr }\widetilde{f}=\widetilde{\psi}$ on $\TT$. Define now $f(z):=\overline{\widetilde{f}\left(\frac{r_0}{\overline{z}}\right)}$ for all $z\in \Omega$. Then, $f\in H^p_{\nu}(\Omega)$, $f$ is continuous on $\overline{\Omega}$ and $\mbox{Re tr }f=\psi$ on $r_0\TT$. \par
\noindent {\bf Step 3: } Assume now that $\Omega=\GG$ is a circular domain, as in \eqref{circ_dom}. Extend $\nu$ to a function $\widetilde{\nu}\in W^{1,\infty}_{\RR}(\CC)$ satisfying the properties of Lemma \ref{extension}. If $\psi\in C^1_{\RR}(\TT)$, step $1$ provides a function $f\in H^p_{\widetilde{\nu}}(\DD)$, continuous on $\overline{\DD}$, and such that $\mbox{Re tr }f=\psi$ on $\TT$. The restriction of $f$ to $\GG$ belongs to $H^p_{\nu}(\GG)$ and satisfies all the requirements. If $\psi\in C^1_{\RR}(a_j+r_j\TT)$, argue similarly using Step $2$ instead of Step $1$. \par
\noindent {\bf Step 4: } Finally, in the general case where $\Omega$ is a Dini-smooth $n$-connected domain, $\Omega$ is conformally equivalent to a circular domain $\GG$, via a confomal map which is $C^1$ up to the boundary of $\Omega$, and we conclude the proof using Step $3$.
\endpf\par
\noindent Let $\Omega_1,\Omega_2$ be domains in $\CC$ and $\phi:\Omega_1\to\Omega_2$ be analytic with $\phi\in W_{\Omega_2}^{1,\infty}(\Omega_1)$. The adjoint of the operator $C_{\phi}$ will play an important role in the following arguments. Note first that, by Proposition \ref{continuouscphi}, $C_{\phi}^{\ast}$ is a bounded linear operator from $(H^p_{\nu\circ\phi}(\Omega_1))^{\prime}$ to $(H^p_{\nu}(\Omega_2))^{\prime}$. Moreover:
\begin{lemma} \label{adjoint}
For all $z\in \Omega_1$, $C_{\phi}^{\ast}({\mathcal E}^{\nu\circ\phi}_z)={\mathcal E}^{\nu}_{\phi(z)}$ and $C_{\phi}^{\ast}({\mathcal F}^{\nu\circ\phi}_z)={\mathcal F}^{\nu}_{\phi(z)}$.
\end{lemma}
\beginpf Let $f\in H_{\nu}^p(\Omega_2)$. Then
$$
\langle C_{\phi}^{\ast}({\mathcal E}^{\nu\circ\phi}_z),f\rangle =\langle {\mathcal E}^{\nu\circ\phi}_z,C_{\phi}(f)\rangle=\langle {\mathcal E}^{\nu\circ\phi}_z, f\circ\phi\rangle=\mbox{ Re }f(\phi(z))=\langle {\mathcal E}^{\nu}_{\phi(z)},f\rangle,
$$
and the argument is analogous for ${\mathcal F}^{\nu}_z$.
\endpf\par
\begin{theorem} \label{caracinvertible}
Assume that $\Omega_1,\Omega_2$ are bounded Dini-smooth domains. Then, the composition operator $C_\phi:H_\nu^p(\Omega_2)\to H_{\nu\circ\phi}^p(\Omega_1)$ is invertible if, and only if, $\phi$ is a bijection from $\Omega_1$ onto $\Omega_2$. 
\end{theorem}

\beginpf Some ideas of this proof are inspired by \cite[Thm 2.1]{SharmaBhanu}. If $\phi$ is invertible, then $C_{\phi^{-1}}=(C_\phi)^{-1}$.\par
\noindent Assume conversely that $C_\phi$ is invertible. Since $C_{\phi}$ is one-to-one with closed range, for all $L\in (H^p_{\nu\circ\phi}(\Omega_1))^{\prime}$, one has
\begin{equation} \label{closedrange}
\left\Vert C_{\phi}^{\ast}L\right\Vert_{(H^p_{\nu}(\Omega_2))^{\prime}}\gtrsim\left\Vert L\right\Vert_{(H^p_{\nu\circ\phi}(\Omega_1))^{\prime}}.
\end{equation}
Let $z_1$, $z_2\in\Omega_1$ be such that $\phi(z_1)=\phi(z_2)$. Then, by Lemma \ref{adjoint},
$$C_\phi^{*}({\mathcal E}^{\nu\circ\phi}_{z_1})={\mathcal E}^{\nu}_{\phi(z_1)}={\mathcal E}^{\nu}_{\phi(z_2)}=C_\phi^{*}({\mathcal E}^{\nu\circ\phi}_{z_2}).$$
\noindent Since $C_\phi^{*}$ is invertible, it follows that ${\mathcal E}^{\nu\circ\phi}_{z_1}={\mathcal E}^{\nu\circ\phi}_{z_2}$. Similarly, ${\mathcal F}^{\nu\circ\phi}_{z_1}={\mathcal F}^{\nu\circ\phi}_{z_2}$, so that $z_1=z_2$ by Lemma \ref{separation}, and $\phi$ is univalent. \par
\noindent Now, suppose that $\phi$ is not surjective. We claim that 
\begin{equation} \label{notempty}
\partial{\phi(\Omega_1)}\cap\Omega_2\neq\varnothing.
\end{equation}
Indeed, since $\phi$ is analytic and not constant in $\Omega_1$, it is an open mapping, so that $\Omega_2=\phi(\Omega_1)\cup(\Omega_2\cap\partial\phi(\Omega_1))\cup(\Omega_2\backslash\overline{\phi(\Omega_1)})$, the union being disjoint. Assume now by contradiction that \eqref{notempty} is false. Then $\Omega_2$ is the union of the two disjoints open sets in $\Omega_2$, $\phi(\Omega_1)$ and $\Omega_2\backslash\overline{\phi(\Omega_1)}$.  One clearly has $\phi(\Omega_1)\neq \varnothing$. The connectedness of $\Omega_2$ therefore yields that $\Omega_2\backslash\overline{\phi(\Omega_1)}=\varnothing$. In other words, 
\begin{equation} \label{omega2subset}
\Omega_2\subset \overline{\phi(\Omega_1)}.
\end{equation}
But since $\phi$ is assumed not to be surjective, there exists $a\in \Omega_2\setminus \phi(\Omega_1)$, and \eqref{omega2subset} shows that $a\in \Omega_2\cap \partial\phi(\Omega_1)$, which gives a contradiction, since we assumed that \eqref{notempty} was false. Finally, \eqref{notempty} is proved. \par

\medskip

\noindent Let $a\in\partial{\phi(\Omega_1)}\cap\Omega_2$ and $(z_n)_{n\in\NN}$ be a sequence of $\Omega_1$ such that

$$\phi(z_n)\stackrel[n\to\infty]{}\longrightarrow a. $$  

\noindent Up to a subsequence, there exists $z\in \overline{\Omega_1}$ such that $z_n \stackrel[n\to\infty]{}\longrightarrow z$. Note that $z\in \partial\Omega_1$, otherwise $\phi(z)=a	$ which is impossible (indeed, since $a\in \partial\phi(\Omega_1)$ and $\phi(\Omega_1)$ is open, thus $a\notin \phi(\Omega_1)$). Write $\partial\Omega_1=\cup_{j=0}^n \Gamma_j$, where the $\Gamma_j$ are pairwise disjoint Jordan curves, so that $z\in \Gamma_m$ for some $m\in \left\{0,...,n\right\}$.\par
\noindent Now, we claim that 

$$\Vert {\mathcal E}_{z_n}^{\nu\circ\phi}\Vert_{(H^p_{\nu\circ\phi}(\Omega_1))^{\prime}}\stackrel[n\to\infty]{}\longrightarrow +\infty.$$
Indeed, by the very definition of the norm in $(H^p_{\nu\circ\phi}(\Omega_1))^{\prime}$, 

 \begin{equation} \label{normeval}
 \Vert {\mathcal E}_{z_n}^{\nu\circ\phi}\Vert_{(H^p_{\nu\circ\phi}(\Omega_1))^{\prime}}=\sup_{g\in H_{\nu\circ\phi}^p(\Omega_1)\atop{\left\Vert{\text{tr}\,g}\right\Vert_p\leq 1}} \vert \Re g(z_n)\vert.
 \end{equation}
 
\noindent For any $k\in\NN$, there is $f_k\in H_{\nu\circ\phi}^p(\Omega_1)$ such that $\vert f_k(z_n)\vert\stackrel[n\to\infty]{}\longrightarrow k$ and $\left\Vert f_k\right\Vert_{H^p_{\nu\circ\phi}(\Omega_1)}\leq 1$. Indeed, let $\psi_k\in C^{1}_{\RR}(\Gamma_m)$ be such that $\vert\psi_k(z)\vert=2k$ and $\Vert\psi_k\Vert_{L^p(\Gamma_m)}\leq\frac 1C$, where $C$ is the implicit constant in Lemma \ref{solsobolev}. It follows from Lemma \ref{solsobolev} that there is $f_k\in H^p_{\nu\circ\phi}(\Omega_1)$, continuous on $\overline{\Omega_1}$, such that $\Re tr(f_k)=\psi_k$ on $\Gamma_m$ and $\left\Vert f_k\right\Vert_{H^p_{\nu\circ\phi}(\Omega_1)}\leq 1$. Observe that, since $f_k$ is continuous in $\overline{\Omega_1}$, $\left\vert \mbox{Re }f_k(z_n)\right\vert\rightarrow \left\vert \psi_k(z)\right\vert$. As a consequence, there is $N_k\in\NN$ such that 
$$
\left\vert \mbox{Re }f_k(z_n)\right\vert\geq k
$$
for all $n\geq N_k$. Therefore, by \eqref{normeval}, 
$$\Vert {\mathcal E}^{\nu\circ\phi}_{z_n}\Vert_{(H^p_{\nu\circ\phi}(\Omega_1))^{\prime}}\geq k,\,\,n\geq N_k. $$
Thus, $\Vert {\mathcal E}^{\nu\circ\phi}_{z_n}\Vert_{(H^p_{\nu\circ\phi}(\Omega_1))^{\prime}}\to\infty$ as $n\to\infty$, as claimed.\par
\noindent Moreover, by Lemma \ref{contomega}, for all $g\in H_{\nu}^p(\Omega_2)$, ${\cal{E}}^{\nu}_{\phi(z_n)}(g)\stackrel[n\to\infty]{}\longrightarrow{\cal{E}}^{\nu}_{a}(g)$ which proves that 

$$\sup_{n\in\NN}\vert{\cal{E}}^{\nu}_{\phi(z_n)}(g)\vert<\infty.$$

\noindent It follows from the Banach-Steinhaus theorem that the $\Vert{\cal{E}}^{\nu}_{\phi(z_n)}\Vert_{(H^p_{\nu}(\Omega_2))^{\prime}}$ are uniformly bounded. Thus, we have that

$$\frac{\Vert C_\phi^{*} ({\mathcal E}^{\nu\circ\phi}_{z_n})\Vert_{(H^p_{\nu}(\Omega_2))^{\prime}}}{\Vert {\mathcal E}^{\nu\circ\phi}_{z_n}\Vert_{(H^p_{\nu\circ\phi}(\Omega_1))^{\prime}}}=\frac{\left\Vert {\mathcal E}^{\nu}_{\phi(z_n)}\right\Vert_{(H^p_{\nu}(\Omega_2))^{\prime}}}{\Vert {\mathcal E}^{\nu\circ\phi}_{z_n}\Vert_{(H^p_{\nu\circ\phi}(\Omega_1))^{\prime}}}\stackrel[n\to\infty]{}\longrightarrow 0,$$

\noindent which contradicts \eqref{closedrange}. We conclude that $\phi$ is surjective. 

\endpf\par

\begin{remark}
\begin{itemize}
\item[$1.$]
To our knowledge, the conclusion of Theorem \ref{caracinvertible} is new, even for Hardy spaces of analytic functions when $\Omega_1$ or $\Omega_2$ are multi-connected.
\item[$2.$]
The proof of Theorem \ref{caracinvertible} does not use the explicit description of the dual of $H_{\nu}^p(\Omega_2)$ and $H_{\nu\circ\phi}^p(\Omega_1)$. Such a description exists when $\Omega_1$ and $\Omega_2$ are simply connected (see \cite[Thm 4.6.1]{BLRR}).
\end{itemize}
\end{remark}

\noindent It follows easily from Theorem \ref{caracinvertible} and Lemma \ref{boundedgpalpha} that:
\begin{corollary} \label{invertiblegpalpha}
Let $\Omega_1,\Omega_2$ be bounded Dini-smooth domains and $\phi\in W^{1,\infty}_{\Omega_2}(\Omega_1)$ be analytic in $\Omega_1$. Let $\alpha\in L^{\infty}(\Omega_2)$. Then $C_\phi: G^p_{\alpha}(\Omega_2)\rightarrow G^p_{\tilde{\alpha}}(\Omega_1)$ is an isomorphism if and only if $\phi$ is a bijection from $\Omega_1$ onto $\Omega_2$.
\end{corollary}

\noindent The characterizations given in Theorem \ref{caracinvertible} and Corollary \ref{invertiblegpalpha} are the same as in the analytic case when $\Omega=\DD$ (see \cite{schwa}). 


\section{Isometries and composition operators on generalized Hardy spaces}\label{sec:isom}

\noindent Throughout this section, $\Omega$ will denote the unit disc $\DD$ or the annulus $\AAA$ and $G^p_{\alpha}(\Omega)$ is equipped with the norm:
$$
\left\Vert \omega\right\Vert_{G^p_{\alpha}(\Omega)}:=\left\Vert \mbox{tr }\omega\right\Vert_{L^p(\partial\Omega)}
$$
(see \eqref{normtrace}).\par
\noindent Let $\widetilde{\alpha}_0:=\widetilde{\alpha}=(\alpha\circ\phi)\overline{\partial\phi}$ and $\widetilde{\alpha}_{n+1}:=(\widetilde{\alpha}_n\circ\phi)\overline{\partial\phi}$, $n\in\NN$. The arguments below rely on the following observation:

\begin{lemma} \label{preparation}
Let $\Omega$ be the unit disc $\DD$ or the annulus $\AAA$, $\phi:\Omega\rightarrow \Omega$ be a function in $W^{1,\infty}_{\Omega}(\Omega)$ analytic in $\Omega$ and $\alpha\in L^{\infty}(\Omega)$. Assume that $C_\phi$ is an isometry from $G^p_{\alpha}(\Omega)$ to $G^p_{\tilde{\alpha}}(\Omega)$. Then $\phi(\partial\Omega)\subset \partial\Omega$.
\end{lemma}
\beginpf Assume by contradiction that the conclusion does not hold, so that there exists $B_0\subset \partial\Omega$ with $m(B_0)>0$ (where $m$ stands for the $1$-dimensional Lebesgue measure) such that $\phi(B_0)\subset \Omega$.\\

\noindent For $\Omega=\AAA$, either $B_0$ is entirely contained in $\TT$ or in $r_0\TT$ or there exists a Borel set $B\subsetneq B_0$ of positive Lebesgue measure such that $B\subset\TT$. For the last case, we still write $B_0$ instead of $B$ and we can assume without loss of generality that $B_0\subset\TT$. Indeed, if $B_0\subset r_0\TT$, it is enough to use the composition with the inversion $\text{Inv}:z\mapsto \frac{r_0}{z}$ since it is easy to check that the composition operator $C_{\text{Inv}}$ is a unitary operator (invertible and isometric) on $G_{\alpha}^p(\Omega)$ using Proposition 3.2 in \cite{bfl}.\par
\noindent The following argument is reminiscent of \cite{bayart}. Let $\phi_1:=\phi$ and $\phi_{n+1}:=\phi\circ \phi_n$ for all integer $n\geq 1$. Note that $\phi_k(B_0)\subset \Omega$ for all $k\geq 1$. For all integer $n\geq 1$, define
$$
B_n:=\left\{z\in \partial\Omega;\ \phi_n(z)\in B_0\right\}.
$$
Observe that the $B_n$ are pairwise disjoint. Indeed, if $z\in B_n\cap B_m\neq \emptyset$ with $n>m$, then
$$
\phi_n(z)\in B_0\mbox{ and }\phi_m(z)\in B_0,
$$
so that
$$
\phi_{n-m}(\phi_m(z))=\phi_n(z)\in B_0\cap \phi_{n-m}(B_0)\subset B_0\cap \Omega=\emptyset,
$$
which is impossible.\par
\noindent Fix a function $F\in H^p(\Omega)$ such that 

\begin{equation} \label{defF}\left\vert \mbox{tr }F\right\vert \begin{cases}= 1 &\mbox{on } B_0 \\ 
\leq1/2 & \mbox{on } \partial\Omega\setminus B_0. \end{cases}
\end{equation}

\noindent We claim that such a function exists. Indeed, if $\Omega$ is the unit disc $\DD$, the outer function $F$ defined as follows

\begin{equation}\label{spec:outer}
F(z)=\exp\left(\frac{1}{2\pi}\int_{0}^{2\pi}\frac{e^{i\theta}+z}{e^{i\theta}-z} \log\vert g(e^{i\theta})\vert d\theta\right),\,\,z\in\DD,
\end{equation} 

\noindent with $g\in L^p(\TT)$ such that $\left\vert g\right\vert=1$ on $B_0$ and $\left\vert g\right\vert=\frac 12$ on $\partial\Omega\setminus B_0$ satisfies the required conditions.\\

\noindent If $\Omega=\AAA$, we consider the function $f\in H^p(\DD)$ defined as in Equation (\ref{spec:outer}) and $g:\AAA\to\CC$ is the restriction of $f$ to $\AAA$. Observe that $g$ is in $H^p(\AAA)$ for each $p$, since $\left\vert g\right\vert^p=\left\vert f\right\vert^p\leq u$, where $u$ is a harmonic function in $\DD$.
Set $M=\max_{\TT_{r_0}} |g|$. Now let $\widetilde{g}_n(z)=z^n g(z)$, then for $z\in\TT$ we have

\[
 |\widetilde{g}_n(z)| = |z^n g(z)| =|g(z)| = \begin{cases} 1 & \text{for $z\in B_0$} \\ \frac{1}{2} & \text{for $z\in \TT\setminus B_0$} \end{cases}.
\]
For $z\in\TT_{r_0}$, we get
\[
 |\widetilde{g}_n(z)| = |z^n g(z)| = |r_{0}^n| \cdot |g(z)| \le r_{0}^n M.
\]
Now, $r_{0}<1$ so pick $N$ large enough that
\[
 r_{0}^N M < 1/2,
\]
and $F=\widetilde{g}_N$ has the requested properties.\\

\noindent Now, for all integer $j\geq 1$, $F^j\in H^p(\Omega)$ and
$$
\lim_{j\rightarrow +\infty}Ê\left\Vert F^j\right\Vert_{H^p(\Omega)}^p=m(B_0).
$$
Moreover, by the maximum principle, since $F$ is not constant in $\Omega$,
\begin{equation} \label{majorF}
\left\vert F(z)\right\vert<1\mbox{ for all }Êz\in \Omega.
\end{equation}
\noindent By Theorem \ref{difffactor} and Theorem \ref{factorgpalpha} in Appendix \ref{factor} below, for all $j\geq 1$, there exists a function $s_j\in C(\overline{\Omega})$ (indeed, $s_j\in W^{1,r}(\Omega)$ for some $r>2$) with $\mbox{ Re }s_j=0$ on $\partial\Omega$ such that
$$
w_j:=e^{s_j}F^j\in G^p_{\alpha}(\Omega) 
\mbox{ and }
\left\Vert s_j\right\Vert_{L^{\infty}(\Omega)}\leq 4\left\Vert \alpha\right\Vert_{L^{\infty}(\Omega)}.
$$
Thus, since $\mbox{Re }ps_j=0$ on $\partial\Omega$,
\begin{equation} \label{wjcalcul1}
\begin{array}{lll}
\displaystyle \left\Vert w_j\right\Vert_{G^p_{\alpha}(\Omega)}^p  & = &\displaystyle \int_{\partial\Omega} \left\vert \mbox{tr }w_j\right\vert^p= 
\displaystyle \int_{\partial\Omega} \left\vert e^{ps_j}\right\vert \left\vert \mbox{tr }F\right\vert^{jp}\\
& =& \displaystyle \int_{\partial\Omega} \left\vert \mbox{tr }F\right\vert^{jp} = 
\left\Vert F^j\right\Vert_{H^p(\Omega)}^p\rightarrow m(B_0).
\end{array}
\end{equation}
Since $C_\phi$ is an isometry from $G^p_{\alpha}(\Omega)$ to $G^p_{\tilde{\alpha}}(\Omega)$, for all integers $n,j\geq 1$,
\begin{equation} \label{isom1}
\left\Vert (C_\phi)^n w_j\right\Vert_{G^p_{\widetilde{\alpha}_n}(\Omega)}^p=\left\Vert w_j\right\Vert_{G^p_{\alpha}(\Omega)}^p.
\end{equation}
But
\begin{equation} \label{isom2}
(C_{\phi})^n w_j=w_j\circ \phi_n.
\end{equation}
For all $z\in B_n$, $\phi_n(z)\in B_0$ so that, for all $j,n\geq 1$,
\begin{equation} \label{wj1}
\left\vert \mbox{tr }w_j\circ \phi_n(z)\right\vert=\left\vert \mbox{tr }F(\phi_n(z))\right\vert^j =1.
\end{equation}
For all $z\in \partial\Omega\setminus B_n$, $\phi_n(z)\in \overline{\Omega}\setminus B_0$, so that
\begin{equation} \label{wj2}
\left\vert w_j\circ \phi_n(z)\right\vert\leq e^{4\left\Vert\alpha\right\Vert_{L^{\infty}(\Omega)}} \left\vert F(\phi_n(z))\right\vert^j
\end{equation}
if $\phi_n(z)\in \Omega$ and
\begin{equation} \label{wj3}
\left\vert \mbox{tr }w_j\circ \phi_n(z)\right\vert\leq e^{4\left\Vert\alpha\right\Vert_{L^{\infty}(\Omega)}} \left\vert \mbox{tr }F(\phi_n(z))\right\vert^j
\end{equation}
if $\phi_n(z)\in \partial\Omega\setminus B_0$. Gathering \eqref{defF}, \eqref{majorF}, \eqref{isom1}, \eqref{isom2}, \eqref{wj1}, \eqref{wj2} and \eqref{wj3}, one obtains, by the dominated convergence theorem,
\begin{equation} \label{wjcalcul2}
\lim_{j\rightarrow +\infty} \left\Vert \omega_j\right\Vert_{G^p_{\alpha}(\Omega)}^p=m(B_n).
\end{equation}
Comparing \eqref{wjcalcul1} and \eqref{wjcalcul2} yields $m(B_n)=m(B_0)$ for all integer $n\geq 1$. Since $m(B_0)>0$ and the $B_n$ are pairwise disjoint, we reach a contradiction. Finally, $\phi(\partial\Omega)\subset \partial\Omega$.
\endpf\par

%

\subsection{The simply connected case}

\noindent We can now state:
\begin{theorem} \label{isomgpalpha}
Let $\phi:\DD\rightarrow \DD$ satisfying (\ref{phi}),
let $\alpha\in L^{\infty}(\DD)$ and the associated $\tilde{\alpha}$ given by (\ref{alpha}). Then the following assertions are equivalent:
\begin{itemize}
\item[$1.$]
$C_\phi$ is an isometry from $G^p_{\alpha}(\DD)$ to $G^p_{\tilde{\alpha}}(\DD)$,
\item[$2.$]
$C_\phi$ is an isometry from $H^p(\DD)$ to $H^p(\DD)$,
\item[$3.$]
$\phi(0)=0$ and $\phi(\TT)\subset \TT$.
\end{itemize}
\end{theorem}
\beginpf
The equivalence between $2$ and $3$ is contained in \cite[Thm 1]{Forelli}. We now prove that $1$ and $2$ are equivalent. Assume first that $C_\phi$ is an isometry from $G^p_{\alpha}(\DD)$ to $G^p_{\tilde{\alpha}}(\DD)$. Let $F\in H^p(\DD)$. By Theorem \ref{difffactor} in Appendix \ref{factor} below, there exists $s\in C(\overline{\DD})$ such that $\mbox{ Re }s=0$ on $\TT$ and $w:=e^sF\in G^p_{\alpha}(\DD)$. Then, since $\phi(\TT)\subset \TT$ by Lemma \ref{preparation}, one obtains
$$
\left\Vert C_{\phi}F\right\Vert_{H^p(\DD)}=\left\Vert C_{\phi}w\right\Vert_{G^p_{\widetilde{\alpha}}(\DD)}=\left\Vert w\right\Vert_{G^p_{\alpha}(\DD)}=\left\Vert F\right\Vert_{H^p(\DD)}.
$$
Assume now that $C_\phi$ is an isometry on $H^p(\DD)$. Then $3$ holds, so that $\phi(\TT)\subset \TT$. Let $w\in G^p_{\alpha}(\DD)$. Pick up $s\in C(\overline{\DD})$ and $F\in H^p(\DD)$ such that $w=e^sF$, with $\mbox{Re }s=0$ on $\TT$. Since $\left\vert e^s\right\vert=1$ on $\TT$ and $\phi(\TT)\subset \TT$,
$$
\left\Vert w\circ\phi\right\Vert_{G^p_{\widetilde{\alpha}}(\DD)}=\left\Vert e^{s\circ\phi}F\circ\phi\right\Vert_{G^p_{\alpha}(\DD)}=\left\Vert F\circ\phi\right\Vert_{H^p(\DD)}=\left\Vert F\right\Vert_{H^p(\DD)}=\left\Vert w\right\Vert_{G^p_{\alpha}(\DD)}.
$$
\endpf\par
\begin{remark}
The conclusion of Theorem \ref{isomgpalpha} shows that $C_{\phi}$ is an isometry on $H^p(\DD)$ if and only if it is an isometry from $G^p_{\alpha}(\DD)$ to $G^p_{\tilde{\alpha}}(\DD)$ for {\it all } functions $\alpha\in L^{\infty}(\DD)$ and  associated $\tilde{\alpha}$ given by (\ref{alpha}).
\end{remark}
\begin{corollary}
Let $\phi:\DD\rightarrow \DD$ be a function satisfying (\ref{phi}),
$\alpha\in L^{\infty}(\DD)$ and the associated $\tilde{\alpha}$ given by (\ref{alpha}). Then $C_{\phi}$ is an isometry from $G^p_{\alpha}(\DD)$ {\bf onto} $G^p_{\tilde{\alpha}}(\DD)$ if and only if there exists $\lambda\in \CC$ with $\left\vert \lambda\right\vert=1$ such that $\phi(z)=\lambda z$ for all $z\in \DD$. 
\end{corollary}
\beginpf
Assume that $C_{\phi}$ is an isometry from $G^p_{\alpha}$ {\it onto} $G^p_{\tilde{\alpha}}$. Then $C_{\phi}$ is an isomorphism from $G^p_{\alpha}$ {\it onto} $G^p_{\tilde{\alpha}}$, and Theorem \ref{invertiblegpalpha} shows that $\phi$ is bijective from $\DD$ to $\DD$. Moreover, Theorem \ref{isomgpalpha} yields $\phi(0)=0$ and $\phi(\TT)\subset \TT$. These conditions on $\phi$ imply that there exists $\lambda\in \CC$ with $\left\vert \lambda\right\vert=1$ such that $\phi(z)=\lambda z$. The converse is obvious.
\endpf\par

\noindent Let us now turn to the isometry property for the composition operator on $H^p_{\nu}(\DD)$. Here, $H^p_{\nu}(\DD)$ also is equipped with the norm:
$$
\left\Vert f\right\Vert_{H^p_{\nu}(\DD)}:=\left\Vert \mbox{tr }f\right\Vert_{L^p(\TT)}.
$$
We prove:
\begin{proposition} \label{isomhpnu}
Let $\phi:\DD\rightarrow \DD$ be a function in $W^{1,\infty}(\DD)$ analytic in $\DD$. If $C_\phi$ is an isometry from $H^p_{\nu}(\DD)$ to $H^p_{\nu\circ\phi}(\DD)$, then $\phi(\TT)\subset \TT$ and $\phi(0)=0$.
\end{proposition}
\beginpf
Assume that $C_\phi$ is an isometry from $H^p_{\nu}(\DD)$ to $H^p_{\nu\circ\phi}(\DD)$. We first claim that there exists $C>0$ such that, for all $w\in G^p_{\alpha}(\DD)$ and all integer $n\geq 1$,
\begin{equation} \label{cphigalpha}
C^{-1}\left\Vert \mbox{tr }w\right\Vert_{L^p(\TT)}\leq \left\Vert \mbox{tr }(w\circ\phi_n)\right\Vert_{L^p(\TT)}\leq C\left\Vert \mbox{tr }w\right\Vert_{L^p(\TT)},
\end{equation}
where, as in the proof of Lemma \ref{preparation}, $\phi_1:=\phi$ and $\phi_{n+1}:=\phi\circ \phi_n$ for all integer $n\geq 1$. Indeed, let $w\in G^p_{\alpha}(\DD)$ and set $f:=\frac{w+\nu\overline{w}}{\sqrt{1-\nu^2}}$. Then, since $\left\Vert \nu\right\Vert_{L^{\infty}(\DD)}<1$, one has, for almost every $z\in \TT$,
\begin{equation} \label{equivfw}
\left\vert \mbox{tr }w(z)\right\vert\sim \left\vert \mbox{tr }f(z)\right\vert.
\end{equation}
As a consequence,
\begin{equation} \label{equifwlp}
\left\Vert \mbox{tr }w\right\Vert_{L^p(\TT)}\sim \left\Vert \mbox{tr }f\right\Vert_{L^p(\TT)}
\end{equation}
and, for all $n\geq 1$,
\begin{equation} \label{equifwlpk}
\left\Vert \mbox{tr }(w\circ \phi_n)\right\Vert_{L^p(\TT)}\sim \left\Vert \mbox{tr }(f\circ \phi_n)\right\Vert_{L^p(\TT)},
\end{equation}
where the implicit constant in \eqref{equifwlpk} does not depend on $n$. Since $C_{\phi}$ is an isometry on $H^p_{\nu}(\DD)$, it follows that, for all $n\geq 1$,
\begin{equation} \label{fisom}
\left\Vert \mbox{tr }(f\circ \phi_n)\right\Vert_{L^p(\TT)}=\left\Vert \mbox{tr }f\right\Vert_{L^p(\TT)},
\end{equation}
and \eqref{equifwlp}, \eqref{equifwlpk} and \eqref{fisom} yield \eqref{cphigalpha}.\par
\noindent Let us now establish that $\phi(\TT)\subset \TT$. Argue by contradiction and let $B_n$ (for all $n\geq 0$) as in the proof of Lemma \ref{preparation}. Consider a function $F\in H^p(\DD)$ and define the functions $s_j$ and $w_j$ as in the proof of Lemma \ref{preparation}. By \eqref{cphigalpha}, for all integers $n,j\geq 1$,
\begin{equation} \label{isom1hpnu}
\left\Vert \mbox{tr } w_j\circ\phi_n\right\Vert_{L^p(\TT)}^p\sim \left\Vert \mbox{tr }w_j\right\Vert_{L^p(\TT)}^p.
\end{equation}
But, as already seen,
$$
\left\Vert \mbox{tr } w_j\circ\phi_n\right\Vert_{L^p(\TT)}^p\rightarrow m(B_n),
$$
so that, by \eqref{isom1hpnu}, $m(B_n)\gtrsim m(B_0)$ for all integer $n\geq 1$. Since $m(B_0)>0$ and the $B_n$ are pairwise disjoint, we reach a contradiction. Finally, $\phi(\TT)\subset \TT$.\par
\noindent Let us now prove that $\phi(0)=0$. Recall now that, since $C_{\phi}$ is an isometry, for all functions $f,g\in H^p_{\nu}(\DD)$, see \cite[Lem. 1.1]{mv}:
\begin{equation} \label{fg}
\int_{\TT} (\mbox{tr }f\circ \phi)\left\vert \mbox{tr }g\circ\phi\right\vert^{p-2}\overline{\mbox{tr }g\circ\phi}=\int_{\TT} \mbox{tr }f\left\vert \mbox{tr }g\right\vert^{p-2}\overline{\mbox{tr }g} \, .
\end{equation}
Applying \eqref{fg} with $g=1$, one obtains, for all $f\in H^p_{\nu}(\DD)$,
\begin{equation} \label{integralf}
\int_{\TT} (\mbox{tr }f\circ \phi)=\int_{\TT} \mbox{tr }f.
\end{equation}
Let $u\in L^p_{\RR}(\TT)$ and $f\in H^p_{\nu}(\DD)$ such that $\mbox{Re tr }f=u$. Taking the real part in the both sides of \eqref{integralf} yields
\begin{equation} \label{integuphi}
\int_{\TT} u\circ \phi=\int_{\TT} u.
\end{equation}
Since this is true for all $u\in L^p_{\RR}(\TT)$, one obtains that \eqref{integralf} holds for all $f\in H^p(\DD)$ (write $\mbox{ tr }f=u+iv$ and apply \eqref{integuphi} with $u$ and $v$), and this yields $\phi(0)=0$ ($f(z)=z$ in (\ref{integralf})).
\endpf\par

\noindent As a corollary of Proposition \ref{isomhpnu}, we characterize isometries from $H^p_{\nu}(\DD)$ {\it onto} $H^p_{\nu\circ\phi}(\DD)$:

\begin{corollary}
Let $\phi:\DD\rightarrow \DD$ be a function in $W^{1,\infty}(\DD)$ analytic in $\DD$. Then $C_{\phi}$ is an isometry from $H^p_{\nu}(\DD)$ {\it onto }$H^p_{\nu\circ\phi}(\DD)$ if and only if there exists $\lambda\in \CC$ with $\left\vert \lambda\right\vert=1$ such that $\phi(z)=\lambda z$ for all $z\in \DD$.
\end{corollary}
\beginpf
Proposition \ref{isomhpnu} shows that $\phi(\TT)\subset \TT$ and $\phi(0)=0$, Theorem \ref{caracinvertible} ensures that $\phi$ is a bijection from $\DD$ onto $\DD$, and the conclusion readily follows. 
\endpf\par
Note that we do not know how to characterize those composition operators which are isometries from $H^p_{\nu}(\DD)$ to $H^p_{\nu\circ\phi}(\DD)$.
\subsection{The case of doubly-connected domains}

In the annular case, we obtain a complete description of 
the composition operators which are isometries on generalized Hardy spaces on $\AAA$. Before stating this result, we check:
\begin{lemma} \label{preparationbis}
Let $\phi:\AAA\rightarrow \AAA$ be analytic with $\phi\in W_{\AAA}^{1,\infty}(\AAA)$. 
\begin{itemize}
\item[$1.$]
If $C_{\phi}$ is an isometry from $G^p_{\alpha}(\AAA)$ into $G^p_{\tilde{\alpha}}(\AAA)$, then $\phi(\partial\AAA)\subset \partial\AAA$.
\item[$2.$]
If $C_{\phi}$ is an isometry from $H^p_{\nu}(\AAA)$ into $H^p_{\tilde{\nu}}(\AAA)$, then $\phi(\partial\AAA)\subset \partial\AAA$.
\end{itemize}
\end{lemma}
\beginpf Item $1$ is already stated in Lemma \ref{preparation}. For item $2$, notice that, if $C_{\phi}$ is an isometry from $H^p_{\nu}(\AAA)$ into $H^p_{\tilde{\nu}}(\AAA)$, then there exists $C>0$ such that, for all $w\in G^p_{\alpha}(\AAA)$,
\begin{equation}\label{ann_cphigalpha}
C^{-1}\left\Vert \mbox{tr }w\right\Vert_{L^p(\partial\AAA)}\leq \left\Vert \mbox{tr }(w\circ\phi_n)\right\Vert_{L^p(\partial\AAA)}\leq C\left\Vert \mbox{tr }w\right\Vert_{L^p(\partial\AAA)}.
\end{equation}
Arguing as in the proof of Proposition \ref{isomhpnu}, one concludes that $\phi(\partial\AAA)\subset\partial\AAA$. 
\endpf\par
\noindent We can now state and prove our description of the composition operators which are isometries on generalized Hardy spaces on $\AAA$:
\begin{theorem} \label{isomannulus}
Let $\phi:\AAA\rightarrow \AAA$ be analytic with $\phi\in W_{\AAA}^{1,\infty}(\AAA)$, $\alpha\in L^{\infty}(\AAA)$ and $\nu$ meeting \eqref{kappa}. The following conditions are equivalent:
\begin{itemize}
\item[$1.$]
$C_{\phi}$ is an isometry from $H^p(\AAA)$ into $H^p(\AAA)$,
\item[$2.$]
$C_{\phi}$ is an isometry from $G^p_{\alpha}(\AAA)$ into $G^p_{\tilde{\alpha}}(\AAA)$,
\item[$3.$]
$C_{\phi}$ is an isometry from $H^p_{\nu}(\AAA)$ into $H^p_{\tilde{\nu}}(\AAA)$,
\item[$4.$]
either there exists $\lambda\in \CC$ of unit modulus such that $\phi(z)=\lambda z$ for all $z\in \AAA$, or there exists $\mu\in \CC$ of unit modulus such that $\phi(z)=\mu\frac{r_0}{z}$ for all $z\in \AAA$.
\end{itemize}
\end{theorem}

\beginpf
The results in Theorem \ref{isom_cns_cont} of Appendix \ref{app} show that $1$ and $4$ are equivalent since $\phi$ is continuous on $\overline{\AAA}$. 
\noindent This shows at once that $1\Rightarrow 2$ and $1\Rightarrow 3$. Assume now that $2$ holds. Then Lemma \ref{preparationbis} shows that $\phi(\partial\AAA)\subset \partial\AAA$. The continuity of $\phi$ on $\partial\AAA$ implies either $\phi(\TT)\subset\TT$ and $\phi(r_0\TT)\subset r_0\TT$ or $\phi(\TT)\subset r_0\TT$ and $\phi(r_0\TT)\subset \TT$. By Remark \ref{ann_contbord} in Appendix \ref{app}, it follows that $4$ holds and therefore item $1$. Finally, if $3$ is true, Lemma \ref{preparationbis} again yields that $\phi(\partial\AAA)\subset \partial\AAA$ continuously and Remark \ref{ann_contbord} entails again item $4$ and thus item $1$ holds.
\endpf\par
Note that, even for Hardy spaces of analytic functions on $\AAA$, the characterization of isometries on $H^p$ given in Theorem \ref{isomannulus} is new.


\section{Compactness of composition operators on Hardy spaces}\label{sec:compact}

Recall that for $\alpha_i\in L^{\infty}(\Omega)$, $\widetilde{\alpha}_i:=(\alpha_i\circ\phi)\overline{\partial\phi}$, $i\in\{1,2\}$. 

\begin{proposition} \label{compact}
Let $\Omega_1,\Omega_2$ be Dini-smooth domains and $\phi:\Omega_1\rightarrow \Omega_2$ be a function in $W_{\Omega_2}^{1,\infty}(\Omega_1)$ analytic in $\Omega_1$. Let $\alpha_1,\alpha_2\in L^{\infty}(\Omega_2)$. Then $C_\phi:G^p_{\alpha_1}(\Omega_2)\rightarrow G^p_{\widetilde{\alpha}_1}(\Omega_1)$ is compact if and only if $C_\phi: G^p_{\alpha_2}(\Omega_2)\rightarrow G^p_{\widetilde{\alpha}_2}(\Omega_1)$ is compact. In particular, if $\alpha\in L^{\infty}(\DD)$, then $C_\phi:G^p_{\alpha}(\DD)\rightarrow G^p_{\widetilde{\alpha}}(\DD)$ is compact if and only if $C_\phi: H^p(\DD)\rightarrow H^p(\DD)$ is compact.
\end{proposition}
We will  use the following observation:
\begin{lemma} \label{convhardy}
Let $\Omega\subset \CC$ be a bounded domain and $\beta\in L^{\infty}(\Omega)$. Let $(w_k)_{k\geq 1}$ be a sequence of locally integrable functions in $\Omega$ and $w\in L^1_{loc}(\Omega)$. Assume that, for all $k\geq 1$, $\overline{\partial}w_k=\beta\overline{w_k}$ and $w_k\rightarrow w$ in $L^1(K)$ for all compact $K\subset \Omega$. Then $\overline{\partial}w=\beta\overline{w}$.
\end{lemma}

\beginpf
Let us first prove Lemma \ref{convhardy}. For all $\varphi\in {\mathcal D}(\Omega)$, since $\varphi$ is compactly supported in $\Omega$,
$$
\displaystyle \langle w,\overline{\partial}\phi\rangle = \displaystyle \lim_{k\rightarrow+\infty} \langle w_k,\overline{\partial}\phi\rangle = - \displaystyle \lim_{k\rightarrow +\infty} \langle \beta \overline{w_k},
\phi\rangle = -\displaystyle \langle \beta \overline{w},\phi\rangle,
$$
which ends the proof.\endpf \par

\medskip

{\sl Proof of Proposition \ref{compact}\,:\,\,} Assume that $C_\phi$ is compact from$G^p_{\alpha_1}(\Omega_2)\rightarrow G^p_{\widetilde{\alpha}_1}(\Omega_1)$. Let $(w^2_k)_{k\geq 1}$ be a bounded sequence in $G^p_{\alpha_2}(\Omega_2)$. By Proposition \ref{alphabeta} in Appendix \ref{factor} below, if $r>2$ is fixed, for all $k\geq 1$, there exists $s_k\in W^{1,r}(\Omega_2)$ with $\mbox{Re }s_k=0$ on $\partial\Omega_2$ and $w^1_k\in G^p_{\alpha_1}(\Omega_2)$ such that
\begin{equation} \label{boundsk}
\left\Vert s_k\right\Vert_{W^{1,r}(\Omega_2)}\lesssim \left\Vert \alpha_1\right\Vert_{L^{\infty}(\Omega_2)}+\left\Vert \alpha_2\right\Vert_{L^{\infty}(\Omega_2)} 
\end{equation}
and
\begin{equation} \label{identity}
w^2_k:=e^{s_k}w^1_k.
\end{equation}
Note that, since for all $k\geq 1$,
$$
\left\Vert s_k\right\Vert_{L^{\infty}(\Omega_2)}\lesssim\left\Vert s_k\right\Vert_{W^{1,r}(\Omega_2)},
$$
the sequence $(w^1_k)_{k\geq 1}$ is bounded in $G^p_{\alpha_1}(\Omega_2)$. Since $C_\phi$ is compact from $G^p_{\alpha_1}(\Omega_2)$ to $G^p_{\widetilde{\alpha}_1}(\Omega_1)$, there exist $\psi:\NN^{\ast}\rightarrow \NN^{\ast}$ increasing and $w^1\in G^p_{\widetilde{\alpha}_1}(\Omega_1)$ such that
$$
w^1_{\psi(k)}\circ \phi\rightarrow w^1\mbox{ in }G^p_{\widetilde{\alpha}_1}(\Omega_1).
$$
By \eqref{boundsk}, up to a second extraction, there exists $s\in W^{1,r}(\Omega_2)$ such that $s_{\psi(k)}\rightarrow s$ weakly in $W^{1,r}(\Omega_2)$ and strongly in $C^{0,\beta}(\overline{\Omega_2})$ where $0<\beta<1-\frac 2r$. Define
$$
w^2:=e^{s\circ\phi}w^1.
$$
For all $k\geq 1$, 
$$
w^2_{\psi(k)}\circ \phi=e^{s_{\psi(k)}\circ\phi}w^1_{\psi(k)}\circ\phi,
$$
and, since $e^{s_{\psi(k)}\circ \phi}\rightarrow e^{s\circ\phi}$ uniformly in $\overline{\Omega_1}$ and $w^1_{\psi(k)}\circ \phi\rightarrow w^1$ in $L^p(\Omega_1)$, one has $w^2_{\psi(k)}\circ\phi\rightarrow w^2$ in $L^p(\Omega_1)$. Hence, Lemma \ref{convhardy} yields that 
$$
\overline{\partial}w^2=\widetilde{\alpha}_2\overline{w^2}.
$$
Pick up a sequence $(\Delta_n)_n$ of domains such that $\overline{\Delta_n}\subset\Omega_1$, $\partial\Delta_n$ is a finite union of rectifiable Jordan curves of uniformly bounded length and each compact subset of $\Omega$
 is eventually contained in $\Delta_n$. For all $n$,
$$
\begin{array}{lll}
\displaystyle \left\Vert w^2\right\Vert_{L^p(\partial\Delta_n)} & \leq  & \displaystyle e^{c\left(\left\Vert \alpha_1\right\Vert_{L^{\infty}(\Omega_2)}+\left\Vert \alpha_2\right\Vert_{L^{\infty}(\Omega_2)}\right)} \left\Vert w^1\right\Vert_{L^p(\partial\Delta_n)}\\
& \leq & \displaystyle e^{c\left(\left\Vert \alpha_1\right\Vert_{L^{\infty}(\Omega_2)}+\left\Vert \alpha_2\right\Vert_{L^{\infty}(\Omega_2)}\right)} \left\Vert w^1\right\Vert_{G^p_{\widetilde{\alpha}_1}(\Omega_1)}.
\end{array}
$$
Thus, $w^2\in G^p_{\widetilde{\alpha}_2}(\Omega_1)$. Moreover, for all $n\in \NN$,
$$
\begin{array}{lll}
\displaystyle \left\Vert w^2_{\psi(k)}\circ \phi-w^2\right\Vert_{L^p(\partial\Delta_n)} & \leq & \displaystyle \left\Vert \left(e^{s_{\psi(k)}\circ \phi}-e^{s\circ\phi}\right)w^1_{\psi(k)}\circ\phi\right\Vert_{L^p(\partial\Delta_n)} + \left\Vert e^{s\circ\phi}\left(w^1_{\psi(k)}\circ\phi-w^1\right)\right\Vert_{L^p(\partial\Delta_n)}\\
& \leq  &  M\left\Vert e^{s_{\psi(k)}\circ \phi}-e^{s\circ\phi}\right\Vert_{L^{\infty}(\Omega_1)}\\
& + & \displaystyle e^{c\left(\left\Vert \alpha_1\right\Vert_{L^{\infty}(\Omega_2)}+\left\Vert \alpha_2\right\Vert_{L^{\infty}(\Omega_2)}\right)} \left\Vert w^1_{\psi(k)}\circ\phi-w^1\right\Vert_{G^p_{\widetilde{\alpha}_1}(\Omega_1)},
\end{array}
$$
where $M:=\sup_{k\geq 1} \left\Vert w^1_{\psi(k)}\circ \phi\right\Vert_{G^p_{\widetilde{\alpha}_1}(\Omega_1)}$. This shows that $w^2_{\psi(k)}\circ \phi\rightarrow w^2$ in $G^p_{\widetilde{\alpha}_2}(\Omega_1)$. Therefore, $C_{\phi}$ is compact from $G^p_{\alpha_2}(\Omega_2)$ to $G^p_{\widetilde{\alpha}_2}(\Omega_1)$. Taking $\alpha_1=\alpha\in L^{\infty}(\DD)$ and $\alpha_2\equiv 0$, we obtain the last characterization. 

\endpf

\noindent The next result is an immediate consequence of Proposition \ref{compact} and Lemma \ref{boundedgpalpha}:
\begin{corollary} \label{compacthpnu}
Under the assumptions of Proposition \ref{compact}, $C_{\phi}$ is compact from $H^p_{\nu}(\DD)$ to $H^p_{\nu\circ \phi}(\DD)$ if and only if $C_{\phi}$ is compact in $H^p(\DD)$.
\end{corollary}

\section{Conclusion}\label{sec:conclu}
We extended to the case of generalized Hardy spaces some well-known properties of  composition operators on classical Hardy spaces of analytic functions $H^p$ for $1 < p < \infty$. Some questions are still open. As mentionned before, it would be interesting to give a complete characterization of isometries among those composition operators on $H^p_{\nu}(\DD)$ spaces. As far as compactness is concerned, we proved that the compactness of $C_{\phi}$ on generalized Hardy spaces is equivalent to the same property on $H^p$. While this property is well-understood in simply-connected domains (\cite[Thm 3.12]{cow},\cite{shapiro1}), multiply-connected situations deserve further investigation. We intend to tackle this issue in a forthcoming work.  
%
%
\appendices

\section{Factorization results for generalized Hardy spaces} \label{factor}

\noindent We first recall here some factorization results relating classical and generalized Hardy spaces $H^p$ and $G_{\alpha}^p$ for $1<p<\infty$.

\begin{proposition} (\cite[Prop. 3.2]{bfl}) \label{easyfactor}

Let $\Omega\subset \CC$ be a bounded Dini-smooth domain, $1<p<+\infty$ and $\alpha\in L^{\infty}(\Omega)$. For all function $w\in G^p_{\alpha}(\Omega)$, there exists $s\in W^{1,r}(\Omega)$ for all $r\in (1,+\infty)$ and $F\in H^p(\Omega)$ such that $w=e^sF$, and $\left\Vert s\right\Vert_{W^{1,r}(\Omega)}\lesssim \left\Vert \alpha\right\Vert_{L^{\infty}(\Omega)}$. Moreover, if $\partial\Omega=\bigcup_{j=0}^n \Gamma_j$ is a finite union of pairwise disjoint Jordan curves, $s$ may be chosen so that, for all $0\leq j\leq n$, $\mbox{Im }s_j=c_j\in \RR$, the $c_j$ are constants, $\sum_{j=0}^n c_j=0$ and one of the $c_j$ can be chosen arbitrarily.
\end{proposition}

\noindent The next Theorem from \cite{bbc} is a kind of converse to Proposition \ref{easyfactor} in the case of Dini-smooth simply connected domain. 

\begin{theorem} (\cite[Thm 1]{bbc}) \label{difffactor}
Let $\Omega\subset \CC$ be a Dini-smooth simply connected domain, $1<p<+\infty$ and $\alpha\in L^{\infty}(\Omega)$. For all function $F\in H^p(\Omega)$, there exists $s\in W^{1,r}(\Omega)$ for all $r\in (1,+\infty)$ such that $w:=e^sF\in G^p_{\alpha}(\Omega)$, $\mbox{Re }s=0$ on $\partial\Omega$ and $\left\Vert s\right\Vert_{W^{1,r}(\Omega)}\lesssim \left\Vert \alpha\right\Vert_{L^{\infty}(\Omega)}$. 
\end{theorem}

\noindent We extend Theorem \ref{difffactor} (Theorem $1$ in \cite{bbc}) to the case of $n$-connected Dini smooth domains. 

\begin{theorem} \label{factorgpalpha}
Let $\Omega\subset \CC$ be a $n$-connected Dini smooth domain. Let $F\in H^p(\Omega)$, $\alpha\in L^{\infty}(\Omega)$. There exists a function $s\in W^{1,r}(\Omega)$ for all $r\in (2,+\infty)$ such that $\mbox{tr Re}\,s=0$ on $\partial\Omega$, $w=e^sF$ and $\left\Vert s\right\Vert_{W^{1,r}(\Omega)}\lesssim \left\Vert \alpha\right\Vert_{L^{\infty}(\Omega)}$.
\end{theorem}
The proof is inspired by the one of \cite[Theorem 1]{bbc}. By conformal invariance, it is enough to deal with the case where $\Omega=\GG$ is a circular domain. We first assume that $\alpha\in W^{1,2}(\GG)\cap L^{\infty}(\GG)$. For all $\varphi\in W^{1,2}_{\RR}(\GG)$, let $G(\varphi)\in W^{1,2}_{0,\RR}(\GG)$ be the unique solution of
$$
\Delta(G(\varphi))=\mbox{Im }\left(\partial(\alpha e^{-2i\varphi})\right).
$$
We claim:
\begin{lemma} \label{regulg}
The operator $G$ is bounded from $W^{1,2}_{\RR}(\GG)$ from $W^{2,2}_{\RR}(\GG)$ and compact from $W^{1,2}_{\RR}(\GG)$ to $W^{1,2}_{\RR}(\GG)$.
\end{lemma}
\beginpf
Let $\varphi\in W^{1,2}_{\RR}(\GG)$. As in \cite{bbc}, $\partial(\alpha e^{-2i\varphi})\in L^2(\GG)$ and $\left\Vert \partial(\alpha e^{-2i\varphi})\right\Vert_{L^2(\GG)}\lesssim \left\Vert \varphi\right\Vert_{W^{1,2}(\GG)}$. It is therefore enough to show that the operator $T$, which, to any function $\psi\in L^2_{\RR}(\GG)$, associates the solution $h\in W^{1,2}_{0,\RR}(\GG)$ of $\Delta \psi=h$ is continuous from $L^2(\GG)$ to $W^{2,2}(\GG)$, which is nothing but the standard $W^{2,2}$ regularity estimate for second order elliptic equations (see \cite[Section 6.3, Theorem 4]{evans} and note that $\GG$ is $C^2$). This shows that $G$ is bounded from $W^{1,2}_{\RR}(\GG)$ from $W^{2,2}_{\RR}(\GG)$, and its compactness on $W^{1,2}_{\RR}(\GG)$ follows then from the Rellich-Kondrachov theorem. \endpf\par

\noindent {\sl Proof of Theorem \ref{factorgpalpha} :\,\,}As in the proof of \cite[Theorem 1]{bbc}, Lemma \ref{regulg} entails that $G$ has a fixed point in $W^{1,2}_{\RR}(\Omega)$, which yields the conclusion of Theorem \ref{factorgpalpha} when $\alpha\in W^{1,2}(\Omega)\cap L^{\infty}(\Omega)$, and a limiting procedure ends the proof.\endpf\par

%



\begin{remark}
Note that, contrary to \cite[Thm 1]{bbc}, we did not investigate uniqueness properties of $s$ in Theorem \ref{factorgpalpha}.
\end{remark}

\noindent Combining Proposition \ref{easyfactor} and Theorem \ref{factorgpalpha}, one easily obtains:

\begin{proposition} \label{alphabeta}
Let $\Omega\subset \CC$ be a bounded Dini-smooth domain, $1<p<+\infty$ and $\alpha_1,\alpha_2\in L^{\infty}(\Omega)$. For all $w_1\in G^p_{\alpha_1}(\Omega)$, there exist $s\in W^{1,r}(\Omega)$ for all $r\in (1,+\infty)$ and $w_2\in G^p_{\alpha_2}(\Omega)$ such that $w_1=e^sw_2$, $\mbox{Re }s=0$ on $\partial\Omega$ and $\left\Vert s\right\Vert_{W^{1,r}(\Omega)}\lesssim \left\Vert \alpha_1\right\Vert_{L^{\infty}(\Omega)}+\left\Vert \alpha_2\right\Vert_{L^{\infty}(\Omega)}$.
\end{proposition}
\beginpf
Let $w_1\in G^p_{\alpha_1}(\Omega)$. Proposition \ref{easyfactor} provides $s_1\in W^{1,r}(\Omega)$ for all $r\in (1,+\infty)$ and $F\in H^p(\Omega)$ such that $w_1=e^{s_1}F$, $\mbox{Re }s_1=0$ on $\partial\Omega$ and $\left\Vert s_1\right\Vert_{W^{1,r}(\Omega)}\lesssim \left\Vert \alpha_1\right\Vert_{L^{\infty}(\Omega)}$. Now, Theorem \ref{factorgpalpha} yields $s_2\in W^{1,r}(\Omega)$ for all $r\in (2,+\infty)$ such that  $w_2:=e^{s_2}F\in G^p_{\alpha_2}(\Omega)$, $\mbox{Re }s_2=0$ on $\partial\Omega$ and $\left\Vert s_2\right\Vert_{W^{1,r}(\Omega)}\lesssim \left\Vert \alpha_2\right\Vert_{L^{\infty}(\Omega)}$. Finally, $s:=s_1-s_2$ and $w_1=e^{s}w_2$ meet all the requirements. 
\endpf

\section{Appendix: composition operators on the analytic Hardy spaces of the annulus} \label{app}

Recall first that, if $\phi:\AAA\rightarrow \AAA$ is analytic, then $\phi$ has a nontangential limit almost everywhere on $\partial\AAA$. In the sequel, define, for all $z\in \overline{\AAA}$,
\begin{equation} \label{defphistar}
\phi^{\ast}(z):=\left\{
\begin{array}{ll}
\phi(z) & \mbox{ if }z\in \AAA,\\
\displaystyle \lim_{r\rightarrow 1} \phi(rz) & \mbox{ if }z\in \TT,\\
\displaystyle \lim_{r\rightarrow r_0} \phi(rz) & \mbox{ if }z\in r_0\TT.
\end{array}
\right.
\end{equation}
For an analytic function $\phi:\AAA\rightarrow \AAA$, we define the Borel set, for $r_0\leq \alpha\leq 1$,

$$\Omega_{\phi,\alpha} =\left\{t\in [0,2\pi], \vert\phi^{\ast}(z)\vert=\alpha, z=r_0e^{it}\right\}\cup\left\{t\in [0,2\pi], \vert\phi^{\ast}(z)\vert=\alpha, z=e^{it}\right\}.$$

\noindent With a slight abuse of notation, write

$$\Omega_{\phi,\alpha}\cap r_0\TT =\left\{t\in [0,2\pi], \vert\phi^{\ast}(z)\vert=\alpha, z=r_0e^{it}\right\}$$

\noindent and by 

$$\Omega_{\phi,\alpha}\cap \TT=\left\{t\in [0,2\pi], \vert\phi^{\ast}(z)\vert=\alpha, z=e^{it}\right\}.$$

\noindent Note that $m\left(\Omega_{\phi,\alpha}\right)\in[0,2]$. Let us state:

\begin{theorem} \label{isometry=>inner}
 Let $\phi:\AAA\to\AAA$ be analytic. If $C_\phi$ is an isometry on $H^p(\AAA)$, then 
 
 $$\phi^{\ast}(\partial\AAA)\subset \partial\AAA \,\,\hbox{and}\,\,m\left(\Omega_{\phi,r_0}\right)=m\left(\Omega_{\phi,1}\right)=1.$$
 
\end{theorem}

\beginpf To prove that $\phi^{\ast}(\partial\AAA)\subset\partial\AAA$, we use arguments analogous to the proof of Lemma \ref{preparation}. Suppose toward a contradiction that there is $B_0\subset\AAA$ such that $\phi(B_0)\subset\AAA$. As in the proof of Lemma \ref{preparation}, we can assume that $B_0\subset\TT$. Let $F\in H^p(\AAA)$ satisfying \eqref{defF} with $\Omega=\AAA$. Define also $\phi^{\ast}_1:=\phi^{\ast}$  (see \eqref{defphistar}) and $\phi^{\ast}_{n+1}:=\phi^{\ast}\circ \phi^{\ast}_n$ for all integer $n\geq 1$. For all integer $n\geq 1$, define
$$
B_n:=\left\{z\in \partial\Omega;\ \phi^{\ast}_n(z)\in B_0\right\}.
$$
For all integer $j\geq 1$, $F^j\in H^p(\AAA)$ and for all $n\in\NN$, we have that

$$
m(B_0)=\lim_{j\rightarrow +\infty}Ê\left\Vert F^j\right\Vert_{H^p(\AAA)}^p=\lim_{j\rightarrow +\infty}Ê\left\Vert C_{\phi^{\ast}_n}F^j\right\Vert_{H^p(\AAA)}^p=m(B_n),
$$
 
\noindent which leads to a contradiction since the $(B_n)_n$ are pairwise disjoint. Moreover, we have that 

$$ \Vert C_\phi(Id)\Vert_{p}^p=\int_{0}^{2\pi}|\phi(r_0e^{it})|^p dt + \int_{0}^{2\pi}\left|\phi\left(e^{it}\right)\right|^p dt=r_{0}^p+1,$$

\noindent which implies that 

\begin{equation}\label{eq1}
 r_{0}^p m\left(\Omega_{\phi,r_0}\right)+m\left(\Omega_{\phi,1}\right)=r_{0}^p+1.
 \end{equation}

\noindent Since $m\left(\Omega_{\phi,r_0}\right)+m\left(\Omega_{\phi,1}\right)=2$, we have \eqref{eq1} if and only if $m\left(\Omega_{\phi,r_0}\right)=m\left(\Omega_{\phi,1}\right)=1$.
\endpf 

\noindent For further use, let us consider the subsurface of the logarithm surface, introduced by Sarason in \cite{Sarason65} and defined by 

$$\hat{A}=\{(r,t) : r_0<r<1\,\,\hbox{and}\,\,t\in\mathbb{R}\} .$$

\noindent Let $\psi$ be the function from $\hat{A}$ onto $\AAA$ defined by $\psi(r,t)=re^{it}$. We now recall some definitions needed in the sequel.

\begin{definition}\cite{Sarason65}
A meromorphic function $F$ defined on $\hA$ is said to be modulus automorphic if and only if, for each $(r,t)\in\hA$,

$$ \vert F(r,t+2\pi)\vert= \vert F(r,t)\vert.$$ 

\end{definition}

\noindent By the maximum modulus principle, there is a constant $\lambda$ of unit modulus such that, for all $(r,t)\in\hA$,

$$ F(r,t+2\pi)=\lambda F(r,t).$$ 

\begin{definition}
Such a $\lambda$ is called a multiplier of $F$ and the unique real number $\gamma\in [0,1)$ such that $\lambda=e^{2i\pi\gamma}$ is called the index of $F$.

\end{definition}

\begin{remark}\label{index0}
If $\phi:\AAA\to\AAA$ is analytic, then, $\hphi :\hA\to\AAA$, defined by $\hphi:=\phi\circ\psi$ is analytic on $\hA$ with index equal to $0$ since for each $(r,t)\in\hA$, $\hphi(r,t+2\pi)=\phi(re^{i(t+2\pi)})=\phi(re^{it})$ (the analycity of $\hphi$ following from the analycity of $\phi$ and $\psi$).
\end{remark}

\noindent We recall the definition of the ${\mathcal H}^p_\gamma(\AAA)$ Hardy spaces introduced by Sarason. 

\begin{definition}\label{Hpalpha}
Let $\gamma\in[0,1)$ and $1\leq p<\infty$. A function $F$ is in ${\mathcal H}^p_\gamma(\AAA)$ if and only if $F$ is holomorphic, modulus automorphic of index $\gamma$ such that 

$$ \sup_{r_0<r<1}\int_{0}^{2\pi}\vert F(r,t)\vert^p dt<\infty.$$

\noindent In particular, ${\mathcal H}^p_0(\AAA)$ is the space of all functions on $\hA$ obtained by lifting functions in $H^p(\AAA)$. A function in ${\mathcal H}^p_{\gamma}(\AAA)$ has non-tangential limits at almost every point of $\partial\hA$.

\end{definition}


\noindent We now give the definitions of singular inner functions and outer functions in $\hA$ (see \cite{Sarason65}, p. 27).

\begin{definition}
\begin{itemize}
\item[$1.$]
Let $\mu$ be a finite real Borel measure on the boundary of $\AAA$ and $u$ be the function on $\AAA$ defined by

$$u(re^{it})=\mu\ast K=\int_{-\infty}^{+\infty}K(r,t-s) d\mu(e^{is})+\int_{-\infty}^{+\infty}K\left(\frac{r_0}{r},t-s\right) d\mu(r_0e^{is}),$$

\noindent where the function $K$ is defined on $\hA$ by

$$K(r,t)=\frac{\frac{1}{q_0}\cos\left(\frac{\pi}{q_0}\ln(r/r_0^{1/2})\right)}{\ch\left(\frac{\pi t}{q_0}\right)-\sin\left(\frac{\pi}{q_0}\ln(r/r_0^{1/2})\right)}$$

\noindent and $q_0=-\ln(r_0)$. Let $U$ be the harmonic function defined on $\hA$ obtained by lifting $u$ and $V$ be a harmonic conjugate of $U$. The modulus automorphic function defined by $F=e^{U+iV}$ is said to be associated with the measure $\mu$.
\item[$2.$]
\noindent Say that $F$ is an outer function (resp. a singular inner function) if and only if $\mu$ is absolutely continuous (resp. singular) with respect to the Lebesgue measure.
\end{itemize}
\end{definition}

\begin{remark}\label{sing.sign}
By Corollary $1$ in \cite{Sarason65}, if $\mu$ is a singular measure on $\partial\AAA$ and $F$ is a bounded singular function associated with $\mu$, then $\mu$ is non-positive.
\end{remark}
\noindent Every function $F\in{\mathcal H}^p_\gamma(\AAA)$ has a Riesz-Nevanlinna factorization : $F=HSF_e$ where $H$ is a Blaschke product in $\hA$, $S$ is a singular inner function and $F_e$ is an outer function. We refer to \cite{Sarason65} for more details about Blaschke product on $\hA$.

\begin{proposition}\label{outer}
Any function $F:\hA\to\AAA$ in ${\mathcal H}^p_{\gamma}(\AAA)$ is an outer function.
\end{proposition}

\noindent The proof relies on the following lemma:

\begin{lemma}\label{lim.singular}

If $\mu$ is a non-negative singular measure, then,

$$\lim_{h\longrightarrow 0}\frac{\mu((\theta-h,\theta+h))}{2h}=+\infty\,\,\hbox{for}\,\,\mu\,\,\hbox{almost all}\,\,\theta.$$

\end{lemma}

This can be proved as in Lemma 5.4 Chapter I of \cite{Garnett}, by means of the covering Lemma 4.4 Chapter I of \cite{Garnett}. Let us now turn to the proof of Proposition \ref{outer}.\\

\beginpf By the Riesz-Nevanlinna factorization, one can write $F=HSF_e$. Since $F$ does not vanish on $\hA$ (recall that $F:\hA\to\AAA$), the Blaschke product $H$ is identically equal to $1$. Since $S=e^{U+iV}$, we have that $\vert S\vert=e^{U}$ where 

\begin{equation} \label{defU}
U(r,t)=\int_{-\infty}^{+\infty}K(r,t-s)d\mu(e^{is})+\int_{-\infty}^{+\infty}K\left(\frac{r_0}{r},t-s\right)d\mu(r_0e^{is}).
\end{equation}

\noindent To prove that $F$ has no singular inner factor, we show that for $t\in\mathbb{R}$ fixed,

$$ U(r,t)\stackrel{}{\xrightarrow[r\to 1]{}-\infty}.$$

\noindent Since the second term of the right-hand side of \eqref{defU} is bounded when $r\rightarrow 1$, it suffices to prove that for $t\in\mathbb{R}$ fixed,

$$\int_{-\infty}^{+\infty} K(r,t-s) d\mu(e^{is})\stackrel{}{\xrightarrow[r\to 1]{}-\infty}.$$

\noindent First, note that, since $q_0=-\ln(r_0)$,
$$
\cos\left(\frac{\pi}{q_0}\ln(r/r_0^{1/2})\right)=\cos\left(\frac{\pi}{q_0}\ln(r)+\pi/2\right)=-\sin\left(\frac{\pi\ln(r)}{q_0}\right)
$$
and
$$
\sin\left(\frac{\pi}{q_0}\ln(r/r_0^{1/2})\right)=\sin\left(\frac{\pi}{q_0}\ln(r)+\pi/2\right)=\cos\left(\pi\frac{\ln(r)}{q_0}\right). 
$$
Therefore, $K$ can be rewritten as follows

\begin{equation} \label{rewriteK}
K(r,t-s)=\frac{1}{q_0}\frac{-\sin\left(\frac{\pi\ln(r)}{q_0}\right)}{\ch\left(\frac{\pi(t-s)}{q_0}\right)-\cos\left(\frac{\pi\ln(r)}{q_0}\right)}.
\end{equation}

\noindent By \cite[Thm 7]{Sarason65}, since $F$ is not identically $0$, $S$ is bounded, and Remark \ref{sing.sign} yields that $\mu$ is nonpositive. Then, by the positivity of $K$ on $\hA$ and \eqref{rewriteK}, we have that 

$$\int_{-\infty}^{+\infty} K(r,t-s) d(-\mu(e^{is}))\geq \int_{0}^{2\pi}\frac{1}{q_0}\frac{-\sin\left(\frac{\pi\ln(r)}{q_0}\right)}{\ch\left(\frac{\pi(t-s)}{q_0}\right)-\cos\left(\frac{\pi\ln(r)}{q_0}\right)}d(-\mu(e^{is})).$$

\noindent For fixed $r$, let $s$ be such that 

\begin{equation} \label{t-s}
\vert t-s\vert<\vert\ln(r)\vert=-\ln(r).
\end{equation} 

It follows at once from \eqref{t-s} that there exist two constants $C_1$ and $C_2$ (only depending on $r_0$) such that 

$$-\sin\left(\frac{\pi\ln(r)}{q_0}\right)\geq C_1\vert\ln(r)\vert,$$

\noindent and,

$$\ch\left(\frac{\pi(t-s)}{q_0}\right)-\cos\left(\frac{\pi\ln(r)}{q_0}\right)\leq C_2\ln(r)^2.$$

\noindent Define $h:=1-r$. By the mean-value theorem, there exists $c\in (0,h)$ such that

\begin{equation} \label{MVT}
h\leq \left\vert \ln(1-h)\right\vert=\frac{h}{1-c}\leq \frac{h}{r_0}.
\end{equation}
\noindent As a consequence, there is a constant $C$ (only depending on $r_0$) such that 

$$K(r,t-s)\geq\frac{C}{2(1-r)}=\frac{C}{2h}.$$

\noindent Observe that, by \eqref{MVT}, if $\left\vert t-s\right\vert\leq h$, one has $\left\vert t-s\right\vert\leq \left\vert \ln r\right\vert$, we have that 

\begin{eqnarray*}
\int_{-\infty}^{+\infty} K(r,t-s) d(-\mu(e^{is}))&\geq& \int_{t-h}^{t+h} \frac{C}{2h}d(-\mu(e^{is}))\\
                                                                            &=&C\frac{\mu((t-h,t+h))}{2h}\stackrel{}{\xrightarrow[h\to 0]{}+\infty.}
\end{eqnarray*}

\noindent Consequently, $U(r,t)\stackrel{}{\xrightarrow[r\to 1]{}-\infty}$ which implies that 

\begin{equation} \label{S0}
\vert S(r,t)\vert\stackrel{}{\xrightarrow[r\to 1]{}0}.
\end{equation}

\noindent Moreover, $F_e$ is bounded on $\hat{A}$. Indeed, for all $(r,t)\in \hat{A}$,
$$
\vert F_e(r,t)\vert=\exp\left(\int_{-\infty}^{\infty}K(r,t-s)\ln\vert F(1,s)\vert ds+\int_{-\infty}^{\infty}K\left(\frac{r_0}{r},t-s\right)\ln\vert F(r_0,s)\vert ds\right),
$$
so that
$$
\begin{array}{lll}
\displaystyle \vert F_e(r,t)\vert & \leq  & \displaystyle \exp\left(\Vert \ln |F|\Vert_{L^{\infty}(\partial\hat{A})}\right)\left(\int_{-\infty}^{+\infty} K(r,t-s)+K\left(\frac{r_0}r,t-s\right) ds\right)\\
& \leq & \displaystyle \exp\left(\Vert \ln |F|\Vert_{L^{\infty}(\partial\hat{A})}\right).
\end{array}
$$

Gathering \eqref{S0} and the fact that $F_e$ is bounded, one concludes that, for $t\in\mathbb{R}$,

$$\vert F(r,t)\vert\stackrel{}{\xrightarrow[r\rightarrow 1]{}0},$$

\noindent which is impossible since $F$ takes its values in $\AAA$. So, we deduce that $F$ has no singular inner factor, and thus, $F$ is an outer function.

\endpf

\noindent Now, one can prove the following proposition.

\begin{proposition}\label{possible-cases}
Let $\phi :\AAA\to\AAA$ be an analytic function. Suppose that $C_{\phi}$ is an isometry. Then, only three situations can happen:

\begin{enumerate}
\item[\quad $1)$] $\vert\phi\vert=r_0$ almost everywhere on $\TT_{r_0}$ and $\vert\phi\vert=1$ almost everywhere on $\TT$; 
\item[\quad $2)$] $\vert\phi\vert=1$ almost everywhere on $\TT_{r_0}$ and $\vert\phi\vert=r_0$ almost everywhere on $\TT$; 
\item[\quad $3)$] $m\left(\Omega_{\phi,r_0}\cap \TT_{r_0}\right)=m\left(\Omega_{\phi,r_0}\cap\TT\right)=\frac{1}{2}$.
\end{enumerate}

\end{proposition}

\beginpf Let $\hphi:\hA\to\AAA$ be the function obtained by lifting the analytic function $\phi$. By Definition \ref{Hpalpha}, $\hphi$ is in ${\mathcal H}^p_0(\AAA)$ and by Proposition \ref{outer}, $\hphi$ is an outer function. By Theorem $6$ in \cite{Sarason65}, the index of $\hphi$ (equal to $0$) is congruent modulo $1$ to 

$$\frac{1}{2\pi q_0}\left(\int_0^{2\pi}\ln \vert\hphi(1,t)\vert dt-\int_0^{2\pi}\ln \vert\hphi(r_0,t)\vert dt \right),$$

\noindent which is equal to 

$$\frac{1}{2\pi q_0}\left(\int_0^{2\pi}\ln \vert\phi(e^{it})\vert dt-\int_0^{2\pi}\ln \vert\phi(r_{0}e^{it})\vert dt \right),$$

\noindent which leads after computations to

$$m\left(\Omega_{\phi,r_0}\cap \TT_{r_0}\right)-m\left(\Omega_{\phi,r_0}\cap\TT\right)\equiv 0\,\,\hbox{modulo $1$}.$$

\noindent Since $C_\phi$ is an isometry, by Theorem \ref{isometry=>inner}, $m\left(\Omega_{\phi,r_0}\right)=1$. Writing $s:=m\left(\Omega_{\phi,r_0}\cap \TT_{r_0}\right)$, one has 

$$ s-(1-s)\equiv 0\,\,\hbox{modulo $1$},$$

\noindent which gives that 

$$2s\equiv 0\,\,\hbox{modulo $1$}.$$

\noindent Since $s\in[0,1]$, either $s=0$, either $s=1$, or $s=\frac{1}{2}$, which gives at once the desired conclusion.

\endpf

\begin{remark}\label{firstcase-second}
Note that the second case described in Proposition \ref{possible-cases} follows from the first one after composition by the inversion $z\mapsto\frac{r_0}{z}$.
\end{remark}

\noindent The next theorem completely describes the composition operators on $H^p(\AAA)$ which are isometries in the cases $1$ and $2$ of Proposition \ref{possible-cases}.

\begin{theorem}\label{isom_cns_cont}
Let $\phi :\AAA\to\AAA$ be an analytic function satisfying Case $1)$ (respectively  Case $2)$) of Proposition \ref{possible-cases}. Then,\\
\noindent $C_\phi$ is an isometry on $H^p(\AAA)$ if and only if $\phi(z)=cz$ (respectively $\phi(z)=c\frac{r_0}{z}$) where $c$ is a constant of unit modulus.

\end{theorem}

The proof relies on the following Theorem (Theorem $6$ in \cite{Sarason65}):

\begin{theorem}
Let $F$ be a modulus automorphic outer function. Then,

\begin{enumerate}\label{outer.leq-geq}
\item[\quad$\bullet$] If $\vert F \vert\leq C_1$ almost everywhere on $\partial\hA$, then, $\vert F \vert\leq C_1$ in $\hA$.
\item[\quad$\bullet$] If $\vert F \vert\geq C_2$ almost everywhere on $\partial\hA$, then, $\vert F \vert\geq C_2$ in $\hA$.
\end{enumerate}

\end{theorem}

\beginpf Let us first assume that $\phi$ satisfies Case $1$ of Proposition \ref{possible-cases}.
%
%
\noindent So, the boundary values of $\phi$ are in $L^p(\partial\AAA)$ and

$$\int_0^{2\pi}\ln \vert\phi(e^{it})\vert dt+\int_0^{2\pi}\ln \vert\phi(r_{0}e^{it})\vert dt>-\infty.$$

\noindent Furthermore, by the proof of Proposition \ref{possible-cases}, we know that 

$$\frac{1}{2\pi q_0}\left(\int_0^{2\pi}\ln \vert\phi(e^{it})\vert dt-\int_0^{2\pi}\ln \vert\phi(r_{0}e^{it})\vert dt \right)\equiv 0\,\, \hbox{modulo $1$}.$$

\noindent Applying Theorem $9$ in \cite{Sarason65}, there exists a unique outer function (up to a multiplicative constant of unit modulus) $F\in {\mathcal H}^p_{0}(\AAA)$ such that the modulus of $F$ on $\partial\hA$ is equal almost everywhere to the modulus of $\phi$ on $\partial\AAA$.\\

\noindent Likewise, it follows from Theorem \ref{outer.leq-geq} (with $C_1=1$ and $C_2=r_0$) that $r_0\leq\vert F \vert\leq 1$ in $\overline{\hat{A}}$ and thus, for $(r,t)\in\hA$, $F(r,t)\in\AAA$.\\

 \noindent Now, if we consider the function $\hphi$ defined on $\hA$ obtained by lifting $\phi$ on $\hA$, Proposition \ref{outer} implies that $\hphi$ is an outer function in ${\mathcal H}^p_0(\AAA)$. So, by uniqueness, $F$ is identically equal to $c\hphi$ where $c$ is a unit modulus constant. But, if $\hat{G}$ denotes the function defined on $\hA$ obtained by lifting $z\mapsto z$, then, $\hat{G}$ is also an outer function in ${\mathcal H}^p_0(\AAA)$ with values in $\AAA$. By uniqueness, it follows that, up to a multiplicative constant, $F$ is identically equal to $\hat{G}$ and thus, $\phi$ is the identity on $\AAA$. Argue similarly when $\phi$ satisfies Case $2$. \endpf

\begin{remark}\label{ann_contbord}
More generally, we proved that any analytic function $\phi:\AAA\to\AAA$ and continuous on $\overline{\AAA}$ such that $\phi(\partial\AAA)\subset\partial\AAA$ is either the identity map on $\AAA$ or the inversion map ($\phi(z)=\frac{r_0}{z}$), up to a unimodular multiplicative constant. 
\end{remark}


\begin{thebibliography}{AAA}
\bibitem{ahlfors} L. Ahlfors, {\it Lectures on quasiconformal
mappings}, Wadsworth and Brooks/Cole Advanced Books and Software,
1987.
\bibitem{ap} K. Astala, L. P\"aiv\"arinta, Calder\'on's inverse
conductivity problem in the plane, {\it Ann. of Math.} {2} (16),
no. 1, 265--299, 2006.
\bibitem{alessandrini-rondi} G. Alessandrini, L. Rondi, Stable determination of a crack in a planar inhomogeneous conductor, {\it SIAM J. Math. Anal.} {30} (2), 326--340, 1998.
\bibitem{abr} S. Axler, P. Bourdon, W. Ramey, {\it Harmonic function theory}, Second edition, Graduate Texts in Mathematics {\bf 137}, Springer-Verlag, New York, 2001.
\bibitem{bbc} L. Baratchart, A. Borichev, S. Chaabi, Pseudo-holomorphic functions at the critical exponent, preprint, http://hal.inria.fr/hal-00824224.
\bibitem{bfl} L. Baratchart, Y. Fischer, J. Leblond, Dirichlet/Neumann problems and Hardy classes for the planar conductivity equation, {\it Compl. Var. Elliptic Eq.}, to appear.
\bibitem{BLRR} L. Baratchart, J. Leblond, S. Rigat and E. Russ, Hardy spaces of the conjugate {B}eltrami equation, {\it J. Funct. Anal.} {259} (2), 384--427, 2010.
\bibitem{bayart} F. Bayart, Similarity to an isometry of a composition operator, {\it Proc. Amer. Math. Soc.} {131} (6), 1789--1791, 2002.
\bibitem{bn} L. Bers, L. Nirenberg, On a representation theorem for
linear elliptic systems with discontinuous coefficients and its
applications, {\it Conv. Int. EDP}, Cremonese, Roma, 111--138, 1954.
\bibitem{SharmaBhanu} U.~Bhanu, S.D. Sharma, Invertible and isometric composition operators on vector-valued {H}ardy spaces, {\em Bull. Korean Math. Soc.}, 41, 413--418, 2004.
\bibitem{Boyd76} D.M. Boyd, Composition operators on {$H^{p}(A)$}, {\em Pacific Journal of Mathematics}, 62 (1), 55--60, 1976.
\bibitem{cow} C.C. Cowen, B.D. MacCluer, Composition operators on spaces of analytic functions, {\em Studies in Advanced Mathematics}, 1995. 
\bibitem{duren} P. L. Duren, {\it Theory of $H^p$ spaces}, Pure and Applied Mathematics, Vol. 38 Academic Press, New York-London, 1970.
\bibitem{EfendievRuss} M. Efendiev, E. Russ,
Hardy spaces for the conjugated Beltrami equation in a doubly connected domain,
{\it J. Math. Anal. and Appl.}, 383,   439--450, 2011.
\bibitem{evans} L.C. Evans, {\it Partial Differential Equations}, Amer. Math. Soc., 1998.
\bibitem{fischer} Y. Fischer, Approximation dans des classes de fonctions analytiques g\'en\'eralis\'ees et r\'esolution de probl\`emes inverses pour les tokamaks,
PhD Thesis, Univ. Nice-Sophia Antipolis, 2011.
\bibitem{flps} Y. Fischer, J. Leblond, J.R. Partington, E. Sincich, 
Bounded extremal problems in Hardy spaces for the conjugate Beltrami equation in simply connected domains, {\it Appl. Comp. Harmo. Anal.}, 31, 264--285, 2011. 
\bibitem{Forelli} F. Forelli,The isometries of {$H^{p}$},{\em Canad. J. Math.}, {16}, 721--728, 1964.
\bibitem{Garnett} J.B. Garnett, {\it Bounded analytic functions}, Pure and Applied Mathematics {96}, 1981.
\bibitem{krav} V.V. Kravchenko, {\it Applied pseudoanalytic function theory}, Frontiers in Math., Birkh\"auser Verlag, 2009. 
\bibitem{llqr} P. Lef\`evre, D. Li, H. Queff\'elec, L. Rodr\'\i guez-Piazza, Compact composition operators on the
Dirichlet space and capacity of sets of
contact points, {\it J. Funct. Anal.} {\bf 264}, 4, 895--919, 2013.
\bibitem{mv} M.J. Martin, M. Vukotic, Isometries of some classical function spaces among the composition operators, in {\it Recent advances in operator-related function theory}, {\it Contemp. Math.} {393}, 133--138, 2006.
\bibitem{musaev} K.M. Musaev, Some classes of generalized analytic
functions, {\it Izv. Acad. Nauk Azerb. S.S.R.}, 2, 40-46, 1971 (in Russian).
\bibitem{nord} E. Nordgren, P. Rosenthal, F.S. Wintrobe, Invertible composition operators on {$H^p$}, {\em J. Funct. Anal.}, 73 (2), 324--344, 1987.  
\bibitem{rudin} W. Rudin, Analytic functions of class $H_p$, {\it Trans. Amer. Math. Soc.} {78}, 46Ð66, 1955. 
\bibitem{Sarason65} D. Sarason, {\it The {$H^{p}$} spaces of an annulus}, Memoirs of the A.M.S. {56}, 1965.
\bibitem{shapiro1} J.H. Shapiro, The essential norm of a composition operator, {\it Ann. of
Math.} {125}, 375--404, 1987.
\bibitem{shapiro2} J. H. Shapiro, {\it Composition Operators and Classical Function Theory}, Springer, New York, 1993.
\bibitem{shapsmith} J.H. Shapiro, W. Smith, Hardy spaces that support no compact composition operators,{\em J. Funct. Anal.} 205 (1), 62-89, 2003. 
\bibitem{schwa} H.J. Schwartz, {\it Composition operators on {$H^p$}}, Thesis, University of Toledo, Toledo, Ohio, 1969.
\bibitem{vekua} I.N. Vekua, {\it Generalized Analytic Functions},
Addison-Wesley, 1962.

\end{thebibliography}
\end{document}